\documentclass{amsart}
\usepackage{amssymb, amsmath, latexsym}

\renewcommand{\baselinestretch}{\baselinestretch}
\renewcommand{\baselinestretch}{1.1}
\author{Constantin-Nicolae Beli}
\title[The $p^n$ torsion of the Brauer group in characteristic $p$]{A
  representation theorem for the $p^n$ torsion of the Brauer group in
  characteristic $p$} 
\date{}

   \def\m{\lim}
   
\def\p{\partial}     
    \def\te{\theta}
   
 \def\({\overline}

\def\){\underline} \def\<{\cdot} \def\go{\mathfrak}

  \def\sbq{\subseteq} \def\spq{\supseteq} 
\def\ti{\times}  \def\oo{{\cal O}} 
   \def\FF{{\mathbb F}}
 \def\ff{\dot{F}} \def\ooo{{\oo^\ti}} 
\def\kk{K^\times} \newcommand\kkk[1]{K^{\times #1}}

 \def\mo{{\rm mod}~}  
  \def\fs{\ff^2}  
\def\p{\go p}    
\def\*{\sharp}  \def\0{} 
 \def\1{^{-1}}  
 \def\[{\prec} \def\]{\succ} 
\def\bmat{\left(\begin{array}} \def\emat{\end{array}\right)}

\def\N{{\rm N}}  
  
 \def\m2{~(\mo 2)} \def\no{\noindent}
 \def\btm{\begin{thm}}
\def\etm{\end{tm}}
 \def\blem{\begin{lem}}
\def\elem{\end{lem}}

\newtheorem{theorem}{Theorem}[section]
\newtheorem{proposition}[theorem]{Proposition}
\newtheorem{lemma}[theorem]{Lemma}
\newtheorem{definition}{Definition}
\newtheorem{corollary}[theorem]{Corollary}

\newtheorem{bof}[theorem]{}
\newtheorem{teorema}{Theorem}

\def\qed{\mbox{$\Box$}\vspace{\baselineskip}}
\def\pf{$Proof.$\,\,} 
\def\bco{\begin{corollary}} \def\eco{\end{corollary}} 
\def\bdf{\begin{definition}} \def\edf{\end{definition}} 
\def\btm{\begin{theorem}} \def\etm{\end{theorem}} 
\def\bpr{\begin{proposition}} \def\epr{\end{proposition}}  
\def\blm{\begin{lemma}} \def\elm{\end{lemma}} 
\def\bff{\begin{bof}\rm} \def\eff{\end{bof}}
\def\btr{\begin{teorema}} \def\etr{\end{teorema}}

\def\de{\newcommand} \de\tm[1]{{\no\bf Theorem~#1}} 
\def\mb{\mathbb} 
 \def\QQ{{\mb Q}}  \def\ZZ{{\mb Z}}
\def\NN{{\mb N}} \def\FF{{\mb F}}

\de\lm[1]{{\no\bf Lemma~#1}}
\de\df[1]{{\no\bf Definition~#1}} \de\co[1]{{\no\bf Corollary~#1}}
\de\tp[1]{\te (#1 )} \de\ts[1]{\te (O^-(#1 ))} \de\ty[1]{\te
(O(#1 ))} \de\tx[1]{\te (#1 )} \de\up[1]{(1+\p^{#1} )\fs}
 \de\ups[1]{((1+\p^{#1})\fs )^*} \de\upo[1]{(1+\p^{#1} )\ooo^2}
\de\upon[2]{(1+\p^{#1})\ooo^2\cap\N (#2 )}
\de\lr[1]{\longrightarrow^{\!\!\!\!\!\!\!\! #1}}
\de\lf[1]{\longleftarrow^{\!\!\!\!\!\!\!\! #1}}
\de\si[1]{\sim^{\!\!\!\!\! #1}} \de\apr[1]{\approx^{\!\!\!\!\! #1}}
\de\leg[2]{\left(\frac {#1}{#2}\right)}
\DeclareMathOperator\Br{Br}

\de\Brr[1]{{}_{#1}\Br}
\newcommand{\diff}{\mathop{}\!\mathrm{d}}

\DeclareMathOperator\Ima{Im}
\DeclareMathOperator\coker{coker}
\DeclareMathOperator\car{char}

\DeclareMathOperator\dlog{dlog}

\begin{document}

\maketitle

\begin{abstract}
If $K$ is a field of characteristic $p$ then the $p$-torsion of the
Brauer group, $\Brr p(K)$, is represented by a quotient of the group
of $1$-forms, $\Omega^1(K)$. Namely, we have a group isomorphism
$$\alpha_p:\Omega^1(K)/\langle\diff a,\, (a^p-a)\dlog b\, :\, a,b\in
K,\, b\neq 0\rangle\to\Brr p(K),$$
given by $a\diff b\mapsto [ab,b)_p$ $\forall a,b\in K$, $b\neq
0$. Here $[\cdot,\cdot )_p$ is the Artin-Schreier symbol.

In this paper we prove that for every $n\geq 1$ we have a
representation of $\Brr{p^n}(K)$ by a quotient of
$\Omega^1(W_n(K))$, where $W_n(K)$ is the truncation of length $n$ of
the ring of $p$-typical Witt vectors, i.e. $W_{\{1,p,\ldots,p^{n-1}\}}(K)$.
\end{abstract}
\section{Introduction}

The Hilbert symbol of degree $p$, $(\cdot,\cdot
)_p:\kk/\kkk p\times\kk/\kkk p\to\Brr p(K)$, which is defined when $\car
K\neq p$ and $\mu_p\sbq K$, has an analogue in characteristic $p$. If
$\wp$ is the Artin-Schreier map, $x\mapsto x^p-x$, then we have the
Artin-Schreier symbol
$$[\cdot,\cdot )_p:K/\wp (K)\times\kk/\kkk p\to\Brr p(K),$$
given for every $a,b\in K$, $b\neq 0$ by $[a,b)_p=[A_{[a,b)_p}]$,
which is the class in the Brauer group of the central simple algebra (c.s.a.)
$A_{[a,b)_p}$ generated by $x,y$ with the relations 
$$\wp (x)=x^p-x=a,\quad y^p=b,\quad yxy^{-1}=x+1.$$

The representation of $\Brr p(K)$ involves the groups of $1$-forms,
$\Omega^1(K)$. Recall, over an arbitrary abelian ring $R$,
$\Omega^1(R)$ is the $R$-module generated by $\diff a$ with $a\in R$,
subject to $\diff (a+b)=\diff a+\diff b$ and $\diff (ab)=a\diff
b+b\diff a$. By [GS, Lemma 9.2.1] we have a group morphism $\gamma
:\Omega^1(K)\to\Omega^1(K)/\diff K$, called the inverse Cartier
operator, or $C^{-1}$, given by $a\diff b\mapsto a^pb^{p-1}\diff
b$. An important property of $\gamma$ is a the existence of the
following exact sequence:
$$1\to\kk\xrightarrow
p\kk\xrightarrow\dlog\Omega^1(K)\xrightarrow{\gamma
-1}\Omega^1(K)/\diff K.$$
This result, [GS, Theorem 9.2.2], is due to Jacobson and Cartier. With
the help of this theorem, Kato was able to prove that there is a group
isomorphism
$$\alpha_p:\coker (\gamma -1)=\Omega^1(K)/(\diff K+(\gamma
-1)\Omega^1(K))\to\Brr p(K),$$
given by $a\diff b\mapsto [ab,b)$ $\forall a,b\in K$, $b\neq 0$. (See
[GS, Theorem 9.2.4] and [GS, Proposition 9.2.5].)

Explicitly, the domain of $\alpha_p$ writes as
$\Omega^1(K)/\langle\diff a,\, (a^pb^{p-1}-a)\diff b\, :\, a,b\in
K\rangle$.

Note that if $b\neq 0$ and we write $a=cb$ then
$a^pb^{p-1}-a=(c^p-c)b^{-1}$ so $(a^pb^{p-1}-a)\diff
b=(c^p-c)b^{-1}\diff b=(c^p-c)\dlog b$. Hence in the formula for the
domain of $\alpha_p$ we can replace $(a^pb^{p-1}-a)\diff b$, with
$a,b\in K$ by $(a^p-a)\dlog b=\wp (a)\dlog b$, with $a,b\in K$, $b\neq
0$. Thus the isomorphism $\alpha_p$ is defined as
$$\alpha_p:\Omega^1(K)/\langle\diff a,\,\wp (a)\dlog b\, :\, a,b\in
K,\, b\neq 0\rangle\to\Brr p(K)$$

More generally, if $n\geq 1$ then we have an analogue of the $p^n$th
Hilbert symbol, which is defined in terms of Witt vectors. 

Let $W(K)$ be the ring of $p$-typical Witt vectors over $K$,
i.e. $W_{\{ 1,p,p^2,\ldots\}}(K)$, and let $W_n(K)$ be its truncation
of lenght $n$, i.e. $W_{\{ 1,p,\ldots,p^{n-1}\}}(K)$. (When $n=0$ by
$W_0(K)$ we mean $W_\emptyset (K)=\{ 0\}$.) We have a truncation
morphism $W(K)\to W_n(K)$ given by $(x_0,x_1,\ldots )\mapsto
(x_0,\ldots,x_{n-1})$. More generally, if $m\geq n$ then we have a
trucation map $W_m(K)\to W_n(K)$.

We denote by $F$ the Frobenius endomorphism on $W(K)$ and on
$W_n(K)$, given by $(x_0,x_1,\ldots )\mapsto (x_0^p,x_1^p,\ldots )$,
and by $V$ the Verschiebung map $(x_0,x_1,\ldots )\mapsto
(0,x_0,x_1,\ldots )$, which is additive. Note that for any $n\geq 0$
we can define $V:W_n(K)\to W_{n+1}(K)$ by $(x_0,\ldots,x_{n-1})\mapsto
(0,x_0,\ldots,x_{n-1})$. However, in many cases we will be concerned
with the truncated version, $V:W_n(K)\to W_n(K)$, given by
$(x_0,\ldots,x_{n-1})\mapsto (0,x_0,\ldots,x_{n-2})$. We have the well
known formulas $(Va)b=V(aFb)$ and $FV=VF=p=(x\mapsto px)$,
i.e. $p(x_0,x_1,\ldots )=(0,x_0^p,x_1^p,\ldots )$. More generally,
$(V^ka)b=V^k(aF^kb)$ and $(V^ka)(V^lb)=V^{k+l}(F^laF^kb)$.

For any $a\in K$ we denote by $[a]$ its Teichm\"uller representative
in $W(K)$ or $W_n(K)$, $[a]=(a,0,0,\ldots )$. The map $a\mapsto
[a]$ is multiplicative, but not additive. The zero and unit elements of
the ring of Witt vectors are $0=[0]=(0,0,\ldots )$ and
$1=[1]=(1,0,0,\ldots )$. If $a=(a_0,a_1,a_2,\ldots )$ is a Witt vector
and $b\in K$ then $a[b]=(a_0b,a_1b^p,a_2b^{p^2},\ldots )$.

The map $V^n$ is zero on $W_n(K)$. Moreover, $V^n(W(K))=\{
(0,\ldots,0,a_n,a_{n+1},\ldots )\, :\, a_i\in K\}$, which is the kernel
of the truncation map $W(K)\to W_n(K)$. Therefore $W_n(K)$ can
also be written as $W(K)/V^n(W(K))$. As we will see later, in some
cases it is more advantageous to regard the truncated Witt vectors as
classes of full Witt vectors, especially when we work with truncations
of different lengths.

The Artin-Schreier map on Witt vectors is $\wp =F-1$, given by
$(x_0,x_1,\ldots )\mapsto (x_0^p,x_1^p,\ldots )-(x_0,x_1,\ldots )$. We
have $\ker (\wp :W_n(K)\to
W_n(K))=W_n(\FF_p)\cong\ZZ/p^n\ZZ$. The isomorphism
$\ZZ/p^n\ZZ\to W_n(\FF_p)$ is given by $\overline a\mapsto a\cdot
1=a(1,0,\ldots,0)$. 

We now define the analogue of $(\cdot,\cdot )_{p^n}$ in characteristic
$p$,
$$[\cdot,\cdot )_{p^n}=[\cdot,\cdot )_{K,p^n}:W_n(K)/\wp
(W_n(K))\times\kk/\kkk{p^n}\to\Brr{p^n}(K),$$
given for any $a=(a_0,\ldots,a_{n-1})\in W_n(K)$ and any $b\in\kk$
by $[a,b)_{p^n}=[A_{[a,b)_{p^n}}]$, where $[A_{[a,b)_{p^n}}]$ is a
c.s.a. over $K$ of degree $p^n$ generated by
$x=(x_0,\ldots,x_{n-1})$ and $y$, such that $x$ has mutually commuting
entries, with the relations
$$\wp (x)=Fx-x=a,\quad y^{p^n}=b,\quad yxy^{-1}=x+1.$$
Here $yxy^{-1}:=(yx_0y^{-1},\ldots,yx_{n-1}y^{-1})$ and $x+1$ is the
sum of Witt vectors $(x_0,\ldots,x_{n-1})+(1,0,\ldots,0)$. 
\smallskip

{\bf Note.} The notation $[\cdot,\cdot )$ is not universally
used. Instead of $[(a_0,\ldots,a_{n-1}),b)$ some authors write
$(b,(a_0,\ldots,a_{n-1})]$ or $(b|a_0,\ldots,a_{n-1}]$.

Also the relation $yxy^{-1}=x+1$ from the definition of
$A_{[a,b)_{p^n}}$ is sometimes replaced by $y^{-1}xy=x+1$. With this
alternative definition $A_{[a,b)_{p^n}}$ becomes
$A_{[a,b)_{p^n}}^{op}$ so $[a,b)_{p^n}$ becomes $-[a,b)_{p^n}$, which
is essentially the same thing.
\smallskip

The symbol $[\cdot,\cdot )_{p^n}$ was introduced by Witt in [W] and is
called the Artin-Schreier-Witt symbol. It is bilinear and
$[a,b)_{p^n}$ depends only on the classes of $a\mod\wp (W_n(K))$ and
$b\mod\kkk{p^n}$, which justifies the set of definition. Then $[\wp
(a),b)_{p^n}=[a,b^{p^n})_{p^n}=0$ $\forall a\in W_n(K)$,
$b\in\kk$. Also note that
$[Fa,b)_{p^n}-[a,b)_{p^n}=[Fa-a,b)_{p^n}=[\wp (a),b)_{p^n}=0$ so
$[Fa,b)_{p^n}=[a,b)_{p^n}$. As a consequence,
$[Va,b)_{p^n}=[FVa,b)_{p^n}=[pa,b)_{p^n}=p[a,b)_{p^n}$.

The symbols $[\cdot,\cdot )_{p^n}$ are related to each other by the
formula $[a,b)_{p^n}=[Va,b)_{p^{n+1}}$ $\forall a\in W_n(K)$. (See [W,
Satz 15].) More generally, if $m\geq n$ then
$[a,b)_{p^n}=[V^{m-n}a,b)_{p^m}$. Explicitly,
$[(a_0,\ldots,a_{n-1}),b)_{p^n}=
[(0,\ldots,0,a_0,\ldots,a_{n-1}),b)_{p^m}$. 

We obtain a map between two directed systems.
$$\begin{array}{ccc}
W_n(K)\times\kk&
\xrightarrow{[\cdot,\cdot )_{p^n}}& \Brr{p^n}(K)\\
\hskip -40pt V^{m-n}\times 1\downarrow&{}&\downarrow\\
W_m(K)\times\kk&
\xrightarrow{[\cdot,\cdot )_{p^m}}& \Brr{p^m}(K)
\end{array}.$$
We take the direct limits. The limit of the directed system
$W_0(K)\xrightarrow VW_1(K)\xrightarrow VW_2(K)\xrightarrow V\cdots$
is $CW(K)$, where $(CW(K),+)$ is the group of Witt
covectors. Recall that the elements of $CW(K)$ write as
$(\ldots,a_{-2},a_{-1},a_0)$, with $a_i\in K$, such that $a_i=0$ for
$i\ll 0$. The canonical morphism $\psi_n:W_n(K)\to CW(K)$ is
given by $(a_0,\ldots,a_{n-1})\mapsto (\ldots
0,0,a_0,\ldots,a_{n-1})$. The Frobenius and Verschiebung maps are
defined on $CW(K)$ by
$F(\ldots,a_2,a_1,a_0)=(\ldots,a_2^p,a_1^p,a_0^p)$ and
$V(\ldots,a_2,a_1,a_0)=(\ldots,a_3,a_2,a_1)$. They are are compatible
with the canonical maps, in the sense that $\psi_n(Fa)=F\psi_n(a)$ and
$\psi_n(Va)=V\psi_n(a)$ $\forall a\in W_n(K)$.

Also $\varinjlim\Brr{p^n}(K)=\Brr{p^\infty}(K):=\bigcup_{n\geq
1}\Brr{p^n}(K)$. So we get a symbol
$$[\cdot,\cdot )_{p^\infty}:CW(K)\times\kk\to\Brr{p^\infty}(K).$$
If $a=(\ldots,a_{-2},a_{-1},a_0)\in CW(K)$ with $a_i=0$ for $i\leq -n$
and $b\in\kk$ then $a=\psi_n((a_{-n+1},\ldots,a_0))$. Therefore
$[a,b)_{p^\infty}=[(a_{-n+1},\ldots,a_0),b)_{p^n}$. Since each
$[\cdot,\cdot )_{p^n}$ is bilinear, $[\cdot,\cdot )_{p^\infty}$ is
bilinear as well.

Let $a\in CW(K)$, $b\in\kk$. We write
$a=\psi_n((a_{-n+1},\ldots,a_0))$ for some $n\geq 1$. Then $\wp
(a)=\psi_n(\wp (a_{-n+1},\ldots,a_0))$ so
$[\wp (a),b)_{p^\infty}=[\wp (a_{-n+1},\ldots,a_0),b)_{p^n}=0$. Also,
if $b\in\kkk{p^\infty}:=\bigcup_{n\geq 1}\kkk{p^n}$ then, in
particular, $b\in\kkk{p^n}$ so $[a,b)_{p^\infty}=[
(a_{-n+1},\ldots,a_0),b)_{p^n}=0$. Hence $[a,b)_{p^\infty}=0$ if
$a\in\wp (CW(K))$ or $b\in\kkk{p^\infty}$. Since $[\cdot,\cdot
)_{p^\infty}$ is bilinear, it follows that $[a,b)_{p^\infty}$ depends
only on $a\mod\wp (CW(K))$ and $b\mod\kkk{p^\infty}$. Hence
$[\cdot,\cdot )_{p^\infty}$ can be defined as
$$[\cdot,\cdot )_{p^\infty}:CW(K)/\wp
(CW(K))\times\kk/\kkk{p^\infty}\to\Brr{p^\infty}(K).$$


We also have $[a,b)_{p^n}=[Va,b)_{p^{n+1}}=p[a,b)_{p^{n+1}}$ $\forall
a\in W(K)$, $b\in\kk$. More generally, if $m\geq n$ then
$p^{m-n}[a,b)_{p^m}=[a,b)_{p^n}$. Explicitly,
$p^{m-n}[(a_0,\ldots,a_{m-1}),b)_{p^m}=[(a_0,\ldots,a_{n-1}),b)_{p^n}$. This
is related to the similar relation for Hilbert symbols, $\frac
mn(a,b)_m=(a,b)_n$ (or $(a,b)_m^{\frac mn}=(a,b)_n$, in the more
familiar multiplicative notation).

There is an alternative definition of $[a,b)_{p^n}$ in terms of
cohomology. 

If $a=(a_0,\ldots,a_{n-1})\in W_n(K)$ then $L=K(\wp^{-1}(a))$ is
called an Artin-Schreier-Witt extension of $K$. Explicitly, there is
$\alpha =(\alpha_0,\ldots,\alpha_{n-1})\in W_n(K_s)$ with $\wp (\alpha
)=a$ and
$L=K(\wp^{-1}(a)):=K(\alpha)=K(\alpha_0,\ldots,\alpha_{n-1})$. If
$\alpha,\beta\in\wp^{-1}(a)$ then $\beta -\alpha\in\ker\wp
=W_n(\FF_p)$. Then the extension $L/K$ is Galois and we have an
injective morphism $\chi_a:{\rm Gal}(L/K)\to W_n(\FF_p)$, given by
$\sigma\mapsto\sigma (\alpha )-\alpha$. The field $L=K(\wp^{-1}(a))$
and the morphism $\chi_a$ depend only on $a\mod\wp (W_n(K))$.

Assume that $a_0\notin\wp (K)$. By [W, Satz 13] $\chi_a$ is an
isomorphism so ${\rm Gal}(L/K)\cong W_n(\FF_p)\cong\ZZ/p^n\ZZ$. If
$x=(x_0,\ldots,x_{n-1})$ is a multivariable, then we have a surjective
morphism $K[x]/(\wp (x)-a)\to L$, given by $x\mapsto\alpha$,
i.e. $x_i\mapsto\alpha_i$ $\forall i$. Note that $\wp (x)=a$ writes as
a system of equations $x_i^p-x_i+P_i(x_0,\ldots,x_{i-1})=0$ $\forall
0\leq i\leq n-1$, for some $P_i\in K[X_0,\ldots,X_{i-1}]$. It follows
that $K[x]/(\wp (x)-a)$ has the basis $x_0^{k_0}\cdots
x_{n-1}^{k_{n-1}}$, with $0\leq k_i\leq p-1$. Then $\dim_KK[x]/(\wp
(x)-a)=p^n=\dim_KL$ so $K[x]/(\wp (x)-a)\cong L$. 

If $a\in W_n(K)$ is arbitrary then let $0\leq k\leq n$ be
maximal such that $(a_0,\ldots,a_{k-1})\in\wp (W_k(K))$. Equivalently,
$k$ is maximal such that in the class of $a\mod\wp (W_n(K))$ there is
a Witt vector with $0$ on the first $k$ positions, i.e. $(a+\wp
(W_n(K)))\cap V^k(W_{n-k}(K))\neq\emptyset$. Let
$a'=(a'_0,\ldots,a'_{n-k-1})\in W_{n-k}(K)$ such that $a\equiv
V^ka'\mod\wp (W_n(K))$. If $a'_0\in\wp (K)$ then $a'\equiv Va''\mod\wp
(W_{n-k}(K))$ so $a\equiv V^ka'\equiv V^{k+1}a''\mod\wp (W_n(K))$ for
some $a''\in W_{n-k-1}(K)$, which contradicts the maximality of
$k$. So $a'_0\notin\wp (K)$. Let
$\alpha'=(\alpha'_0,\ldots,\alpha'_{n-k-1})\in W_{n-k}(K_s)$ with $\wp
(\alpha')=a'$. Then $\wp (V^k\alpha')=V^k\wp (\alpha')=V^ka'$. Since
$a\equiv V^ka'\mod\wp (W_n(K))$ we have
$L=K(\wp^{-1}(a))=K(\wp^{-1}(V^ka'))=
K(V^k\alpha')=K(\alpha')=K(\wp^{-1}(a'))$. Since $a'_0\notin\wp (K)$
we have that $\chi_{a'}:{\rm Gal}(L/K)\to W_{n-k}(\FF_p)$ is an
isomorphism. Hence ${\rm Gal}(L/K)\cong
W_{n-k}(\FF_p)\cong\ZZ/p^{n-k}\ZZ$. If $\sigma\in{\rm Gal}(L/K)$ then
$\chi_{a'}(\sigma )=\sigma (\alpha')-\alpha'$ so $\chi_a(\sigma
)=\chi_{V^ka'}(\sigma )=\sigma
(V^k\alpha')-V^k\alpha'=V^k\chi_{a'}(\sigma )$. If we identify
$W_n(\FF_p)$ and $W_{n-k}(\FF_p)$ with $\ZZ/p^n\ZZ$ and
$\ZZ/p^{n-k}\ZZ$ then $W_{n-k}(\FF_p)\xrightarrow{V^k}W_n(\FF_p)$
identifies with $\ZZ/p^{n-k}\ZZ\xrightarrow{p^k}\ZZ/p^n\ZZ$ so we have
$\chi_a(\sigma )=p^k\chi_{a'}(\sigma )$. 

Given a finite Galois extension $L/K$ and $\chi :{\rm
Gal}(L/K)\to\ZZ/p^n\ZZ$ a morphism, we denote by $\tilde\chi\in{\rm
Hom_{cont}}({\rm Gal}(K_s/K),\ZZ/p^n\ZZ )=H^1(K,\ZZ/p^n\ZZ)$ the
induced morphism, given by $\tilde\chi (\sigma )=\chi
(\sigma_{|L})$. By [GS, Remark 4.3.13 2.], we have an isomorphism
$W_n(K)/\wp (W_n(K))\cong H^1(K,\ZZ/p^n\ZZ )$, given by the coboundary
morphism $W_n(K)\to H^1(K,W_n(\FF_p))=H^1(K,\ZZ/p^n\ZZ )$, which comes
from the exact sequence $0\to W_n(\FF_p)\to W_n(K_s)\xrightarrow\wp
W_n(K_s)\to 0$. Explicitly, if $a\in W_n(K)$ and $\alpha\in W_n(K_s)$
such that $\wp (\alpha )=a$ then the element of $H^1(K,W_n(\FF_p))$
corresponding to $a$ is a given by $\sigma\mapsto\sigma (\alpha
)-\alpha$ so it coincides with $\tilde\chi_a$. So the isomorphism
$W_n(K)/\wp (W_n(K))\cong H^1(K,\ZZ/p^n\ZZ )$ is given by
$a\mapsto\tilde\chi_a$.

From the cup product $\cup :H^2(K,\ZZ )\otimes
H^0(K,\kk_s)\to H^2(K,\kk_s)=\Br (K)$ and the coboundary morphism
$\delta :H^1(K,\ZZ/p^n\ZZ )\to H^2(K,\ZZ )$ coming from the exact
sequence $0\mapsto\ZZ\xrightarrow{p^n}\ZZ\to\ZZ/p^n\ZZ\to 0$ we get a
linear map
$$j_n:H^1(K,\ZZ/p^n\ZZ )\otimes K^\times\to\Brr{p^n}(K),\quad
j_n(\psi\otimes b)=\delta (\psi )\cup b.$$

\bpr With the above notations we have

(i) $[a,b)_{p^n}=j_n(\tilde\chi_a\otimes b)=\delta (\tilde\chi_a)\cup
b$.

(ii) $[a,b)_{p^n}=0$ if and only if $b\in{\rm N}_{L/K}(L^\times )$,
where $L=K(\wp^{-1}(a))$.
\epr
\pf Assume first that $a_0\notin\wp (K)$ so that $\chi_a:{\rm
Gal}(L/K)\to\ZZ/p^n\ZZ$ is an isomorphism and we can apply [GS,
Proposition 4.7.3 and Corollary 4.7.5]. By [GS, Proposition 4.7.3] we
have $\delta (\tilde\chi_a)\cup b=[(\chi_a,b)]$, where $(\chi_a,b)$ is
the c.s.a. described in [GS, Proposition 2.5.2]. Namely, $(\chi_a,b)$
is the $K$-algebra generated by $L$ and $y$, subject to the relations
$y^{p^n}=b$ and $y\lambda y^{-1}=\sigma (\lambda )$ $\forall\lambda\in
L$. Here $\sigma$ is the preimage of $1$ under $\chi_a:{\rm
Gal}(L/K)\to W_n(\FF_p)\cong\ZZ/\p^n\ZZ$, i.e. $\sigma$ is given by
$\sigma (\alpha )-\alpha =1$, i.e. by $\alpha\mapsto\alpha +1$, where
$\alpha =(\alpha_0,\ldots,\alpha_{n-1})\in\wp^{-1}(a)$. The relation
$y\lambda y^{-1}=\sigma (\lambda )$ only needs to be verified by the
generators $\alpha =(\alpha_0,\ldots,\alpha_{n-1})$ of $L/K$ so it is
equivalent to $y\alpha y^{-1}=\sigma (\alpha )=\alpha +1$. But, if
$x=(x_0,\ldots,x_{n-1})$ is a multivariable, then we have the
isomorphism $L=K[\alpha]\cong K[x]/(\wp (x)-a)$, given by
$\alpha\mapsto x$, i.e. $\alpha_i\mapsto x_i$ $\forall i$. Hence
$(\chi_a,b)$ is the algebra generated by $x=(x_0,\ldots,x_{n-1})$ and
$y$, where $x_i$'s commute with each other, $\wp (x)=a$, $y^{p^n}=b$
and $yxy^{-1}=x+1$. Hence $(\chi_a,b)=A_{[a,b)_{p^n}}$, so
$\delta(\tilde\chi_a)\cup b=[(\chi_a,b)]=[a,b)_{p^n}$. By [GS,
Corollary 4.7.5] we also have that
$[a,b)_{p^n}=\delta(\tilde\chi_a)\cup b=0$ iff $b\in{\rm
N}_{L/K}(L^\times )$. 

If $a\in W_n(K)$ is arbitrary then let $0\leq k\leq n$ be maximal
such that $(a_0,\ldots,a_{k-1})\in\wp (W_k(K))$. Then there is
$a'=(a'_0,\ldots,a'_{n-k-1})\in W_{n-k}(K)$ with $a'_0\notin\wp (K)$
such that $a\equiv V^ka'\mod\wp (W_n(K))$. We have
$L=K(\wp^{-1}(a))=K(\wp^{-1}(a'))$ and $\chi_a=p^k\chi_{a'}$. Since
$a'_0\notin\wp (K)$ we have
$[a',b)_{p^{n-k}}=\delta'(\tilde\chi_{a'})\cup b$, where
$\delta':H^1(K,\ZZ/p^{n-k}\ZZ )\to H^2(K,\ZZ )$ is the coboundary
morphism obtained from the exact sequence
$0\to\ZZ\xrightarrow{p^{n-k}}\ZZ\to\ZZ/p^{n-k}\ZZ\to 0$, and
$[a',b)_{p^{n-k}}=0$ iff $b\in{\rm N}_{L/K}(L^\times )$. But
$\chi_a=p^k\chi_{a'}$, so $\tilde\chi_a=p^k\tilde\chi_{a'}$, which, by
straightforward calculations, implies that
$\delta(\tilde\chi_a)=\delta'(\tilde\chi_{a'})$, and $a\equiv
V^ka'\mod\wp (W_n(K))$ so $[a,b)_{p^n}=[V^ka',b)_{p^n}=[a',b)_{p^{n-k}}$.
Hence $[a,b)_{p^n}=\delta(\tilde\chi_a)\cup b$ and $[a,b)_{p^n}=0$ iff
$b\in{\rm N}_{L/K}(L^\times )$. \qed

\bpr The group $\Brr{p^n}(K)$ is generated by the image of
$[\cdot,\cdot )_{p^n}$.
\epr
\pf By [GS. Theorem 9.1.4] the map $j_n$ is surjective so
$\Brr{p^n}(K)$ is generated by $j_n(\psi\otimes b)=\delta (\psi )\cup
b$, with $\psi\in H^1(K,\ZZ/p^n\ZZ)$, $b\in\kk$. But every $\psi\in
H^1(K,\ZZ/p^n\ZZ )$ writes as $\tilde\chi_a$ for some $a\in W_n(K)$
and $\delta (\tilde\chi_a)\cup b=[a,b)_{p^n}$. Hence $\Brr{p^n}(K)$ is
generated by $[a,b)_{p^n}$ with $a\in W_n(K)$ and $b\in\kk$. \qed

For the purpose of this paper we only need Proposition 1.1(ii). This
result is very likely already known. However we didn't find it stated
explicitly in the general case in the literature. So we provided a
proof here.
\medskip

If $K$ is a local field then we have another definition of
$[\cdot,\cdot )_{p^n}$, in terms of the local Artin map, with values
in $W_n(\FF_p)$. Namely, if $a\in W_n(K)$ and $b\in\kk$ then we
take $\alpha\in W_n(K_s)$ with $\wp (\alpha )=a$ and we define
$[a,b)_{p^n}:=(b,K(\alpha )/K)(\alpha )-\alpha\in W_n(\FF_p)$. This
new definition of $[\cdot,\cdot )_{p^n}$, with values in $W_n(\FF_p)$,
is related to the initial one, with values in $\Brr{p^n}(K)$, via the
local invariant $inv:\Br (K)\xrightarrow\sim\QQ/\ZZ$. It sends
$\Brr{p^n}(K)$ to $\frac 1{p^n}\ZZ/\ZZ$ so we have
$\Brr{p^n}(K)\cong\frac 1{p^n}\ZZ/\ZZ\cong\ZZ/p^n\ZZ\cong
W_n(\FF_p)$. (See [FV, (7.3)], where we have a general statement for
all cyclic algebras.)


\section{A key lemma}

In this section we prove that for every $b\in\kk$ we have
$[[b],b)_{p^n}=0$. This result, together with [GS, Theorem 9.2.4] and
the basic properties of the $[\cdot,\cdot )_{p^n}$, the bilinearity
and the relations $[\wp (a),b)_{p^n}=0$ (or, equivalently,
$[Fa,b)_{p^n}=[a,b)_{p^n}$) and $[a,b)_{p^n}=[Va,b)_{p^{n+1}}$, are
all the ingredients we will use in this paper.

\blm If $R$ is a ring and ${\mathfrak a}\sbq R$ is an ideal then for
every $n\geq 1$ we have:

(i) $W_n({\mathfrak a})$ is an ideal of $W_n(R)$.

(ii) If $\alpha_h$, with $h\in S$, generate $({\mathfrak a},+)$ then
$V^i[\alpha_h]$, with $h\in S$ and $0\leq i\leq n-1$, generate
$(W_n({\mathfrak a}),+)$.
\elm
\pf (i) We have $W_n({\mathfrak a})=\ker(W_n(R)\to
W_n(R/{\mathfrak a}))$.

(ii) We use the induction on $n$. When $n=1$ we have $W_1({\mathfrak
a})={\mathfrak a}$, which is generated by $[\alpha_h]=\alpha_h$ with
$k\in S$. 

Suppose now that $n>1$ and let $a=(a_0,\ldots,a_{n-1})\in
W_n({\mathfrak a})$. Then $a_0\in{\mathfrak a}$ writes as
$\sum_{h\in S}m_h\alpha_h$ for some $m_h\in\ZZ$ with $m_h=0$ for
almost all $h\in S$. Then $a-\sum_{h\in S}m_h[\alpha_h]$ belongs to
$(W_n({\mathfrak a}),+)$ and its first entry is $a_0-\sum_{h\in
S}m_h\alpha_h=0$. It follows that $a-\sum_{h\in S}m_h[\alpha_h]=Vb$
for some $b=(b_0,\ldots,b_{n-2})\in W_{n-1}({\mathfrak a})$. By the
induction hypothesis, $b$ writes as a linear combination with
coefficients in $\ZZ$ of $V^i[\alpha_h]$, with $h\in S$ and $0\leq
i\leq n-2$. It follows that $Vb$ writes as a linear combination of
$V^i[\alpha_h]$, with $h\in S$ and $1\leq i\leq n-1$. From here we
conclude that $a=\sum_{h\in S}m_k[\alpha_h]+Vb$ writes is a linear
combination of $V^i[\alpha_h]$, with $h\in S$ and $0\leq i\leq
n-1$. \qed

\blm (i) If $n\geq 1$ then $[[b],b)_{p^n}=0$ for every $b\in\kk$.

(ii) More generally, if $k\sbq K$ is a perfect field and $a\in
W_n(bk[b])$ then $[a,b)_{p^n}=0$.
\elm
\pf First we prove (i) at $n=1$. By Proposition 1.1(ii) $[b,b)_p=0$
iff $b\in{\rm N}_{L/K}(L^\times )$, where $L=K(\wp^{-1}(b))$. If
$b\in\wp (K)$ then $[b,b)_p=0$ and $L=K$ so our claim is trivial. If
$b\notin\wp (K)$ then $L/K$ is an Artin-Schreier extension, of degree
$p$. We have $L=K(\alpha )$ for some $\alpha$ with $\wp
(\alpha)=b$. Then the minimal polynomial of $\alpha$ is $X^p-X-b=0$ so
${\rm N}_{L/K}(\alpha )=(-1)^p(-b)=b$. (Including the case when $\car
K=p=2$.) Thus $b\in{\rm N}_{L/K}(L^\times )$.

Next, we prove (i) and (ii) by induction on $n$ in two steps.

{\em Step 1.} We prove that if $n\geq 1$ then (i) at $n$ implies (ii)
at $n$. We use Lemma 2.1(ii). Since the ideal ${\mathfrak a}=bk[b]$ of
the ring $R=k[b]$ is generated, as a group, by $cb^h$, with $h\geq 1$
and $c\in k\setminus\{ 0\}$, we get that $W_n(bk[b])$ is generated
by $V^i[cb^h]$, with $h\geq 1$, $c\in k\setminus\{ 0\}$ and $0\leq
i\leq n-1$. Hence $a$ writes as a linear combination with coefficients
in $\ZZ$ of the generators $V^i(cb^h)$ and, by the linearity in the
first variable of $[\cdot,\cdot )_{p^n}$, $[a,b)_{p^n}$ writes as a
linear combination of $[V^i[cb^h],b)_{p^n}$. So it is enough to prove
that $[V^i[cb^h],b)_{p^n}=0$ for $h\geq 1$, $c\in k\setminus\{ 0\}$
and $0\leq i\leq n-1$. But $[V^i[cb^h],b)_{p^n}=p^i[[cb^h],b)_{p^n}$
so we only have to prove that $[[cb^h],b)_{p^n}=0$ (i.e. the case
$i=0$). We write $h=p^sl$ with $s\geq 0$ and $(p,l)=1$. Since $k$ is
perfect we have $c=d^{p^{n+s}}$ for some $d\in k\setminus\{ 0\}$. Then
$[[cb^h],b)_{p^n}=[[d^{p^{n+s}}b^{p^sl}],b)_{p^n}=
[F^s[d^{p^n}b^l],b)_{p^n}=[[d^{p^n}b^l],b)_{p^n}$. By (i) we have
$0=[[d^{p^n}b^l],d^{p^n}b^l)_{p^n}=
p^n[[d^{p^n}b^l],d)_{p^n}+l[[d^{p^n}b^l],b)_{p^n}=
l[[d^{p^n}b^l],b)_{p^n}$. Since also $p^n[[d^{p^n}b^l],b)_{p^n}=0$ and
$(p^n,l)=0$ we get $[[d^{p^n}b^l],b)_{p^n}=0$, as claimed. 

{\em Step 2.} We prove that if $n>1$ then (ii) at $n-1$ implies (i) at
$n$. Let $\alpha =(\alpha_0,\ldots,\alpha_{n-1})\in W_n(K_s)$ with
$\wp (\alpha )=[b]$. By  Proposition 1.1(ii), we must prove that
$b\in{\rm N}_{L/K}(L^\times )$, where $L=K(\wp^{-1}([b]))=K(\alpha )$.
When we identify the first coordinate in the equality $F\alpha -\alpha 
=\wp (\alpha )=[b]$ we get $\alpha_0^p-\alpha_0=b$. We have $\alpha
=(\alpha_0,0,\ldots,0)+(0,\alpha_1,\ldots,\alpha_{n-1})=
[\alpha_0]+V\alpha'$, where $\alpha'=(\alpha_1,\ldots,\alpha_{n-1})\in
W_{n-1}(K_s)$. Then $[b]=\wp (\alpha )=\wp ([\alpha_0])+\wp (V\alpha')$
so $V\wp (\alpha')=\wp (V\alpha')=[b]-\wp
([\alpha_0])=[b]-[\alpha_0^p]+[\alpha_0]$. But
$[b]=[\alpha_0^p-\alpha_0]$, $[\alpha_0^p]$ and $[\alpha_0]$ belong to
$W_n(\alpha_0\FF_p[\alpha_0])$, which, by Lemma 2.1(i), is an ideal
of $W_n(\FF_p[\alpha_0])$. It follows that $V\wp
(\alpha')=[b]-[\alpha_0^p]+[\alpha_0]\in
W_n(\alpha_0\FF_p[\alpha_0])$ so $a:=\wp (\alpha')$ has all the
coordinates in $\alpha_0\FF_p[\alpha_0]$, i.e. $a\in
W_{n-1}(\alpha_0\FF_p[\alpha_0])$. Let now $K'=K[\alpha_0]$. We have
$\wp (\alpha_0)=b$ so $K'=K(\wp^{-1}(b))$. Also
$L=K(\alpha_0,\ldots,\alpha_{n-1})=
K'(\alpha_1,\ldots,\alpha_{n-1})=K'(\alpha')$. Since $\wp
(\alpha')=a\in W_{n-1}(\FF_p[\alpha_0])\sbq W_{n-1}(K')$ we have
$L=K'(\wp^{-1}(a))$. Now $\alpha_0\in K'^\times$, $\FF_p\sbq K'$ is a
perfect field and $a\in W_{n-1}(\alpha_0\FF_p[\alpha_0])$. By (ii) at
$n-1$, this implies that $[a,\alpha_0)_{K',p^{n-1}}=0$. Since
$L=K'(\wp^{-1}(a))$, by Proposition 1.1(ii) we have $\alpha_0\in{\rm
N}_{L/K'}(L^\times )$. There are two cases: 

a) If $b\notin\wp (K)$ then, as seen from the proof of (i) in the case
$n=1$, $\wp (\alpha_0)=b$ implies that $K'=K(\alpha_0)$ is an
Artin-Schreier extension of $K$ and ${\rm N}_{K'/K}(\alpha_0)=b$. But
we also have $\alpha_0\in{\rm N}_{L/K'}(L^\times )$ so $\alpha_0={\rm
N}_{L/K'}(\gamma )$ for some $\gamma\in L^\times$. It follows that
$b={\rm N}_{K'/K}({\rm N}_{L/K'}(\gamma ))= {\rm N}_{L/K}(\gamma )$ so
$b\in{\rm N}_{L/K}(L^\times )$. 

b) If $b\in\wp (K)$ then $K'=K(\wp^{-1}(b))=K$. Hence $\alpha_0\in{\rm
N}_{L/K}(L^\times )$. By the same reasoning, for any other
$\beta=(\beta_0,\ldots,\beta_{n-1})$ with $\wp (\beta )=[b]$ we have
that $\beta_0$ is in the norm group of $L=K(\wp^{-1}([b]))/K$. Let
$0\leq h\leq p-1$. We have $[h]=(h,0,\ldots,0)\in
W_n(\FF_p)=\ker\wp$. So if we take $\beta =\alpha +[h]$ then $\wp
(\beta )=\wp (\alpha )+\wp ([h])=[b]+0=[b]$. By identifying the first
coordinate in the equality $\beta =\alpha +[h]$, we get
$\beta_0=\alpha_0+h$. Hence $\alpha_0+h=\beta_0\in{\rm N}_{L/K}(L^\times
)$. Since $\alpha_0+h\in{\rm N}_{L/K}(L^\times )$ for $0\leq h\leq p-1$,
we have $b=\alpha_0^p-\alpha_0=\alpha_0(\alpha_0+1)\cdots
(\alpha_0+p-1)\in{\rm N}_{L/K}(K^\times )$ and we are done. \qed 

In this paper we only use Lemma 2.2(i). We stated the stronger
result from (ii) only because it makes the induction possible. In
fact, for the induction to work we only need the statement (ii) for
$k=\FF_p$.

Note that the proof of $[[b],b)_{p^n}=0$ was done with rather
rudimentary methods. There is an alternative proof using class field
theory. We can also prove Lemma (ii), but only when $k$ is finite. So
we have $b\in\kk$, $k=\FF_q\sbq K$ for some $p$-power $q$, $a\in
W_n(b\FF_q[b])$ and we want to prove that $[a,b)_{p^n}=0$.

First note that we can reduce to the case when $K=\FF_q(b)$. Indeed,
we have $a\in W_n(\FF_q(b))$ and $\FF_q(b)\sbq K$ so if
$[a,b)_{\FF_q(b),p^n}=0$ then $[a,b)_{K,p^n}=0$, as well. (We have
$A_{[a,b)_{K,p^n}}=A_{[a,b)_{\FF_q(b),p^n}}\otimes_{\FF_q(b)}K$ so if
$A_{[a,b)_{\FF_q(b),p^n}}\cong M_{p^n}(\FF_p(b))$ then also
$A_{[a,b)_{K,p^n}}\cong M_{p^n}(K)$.)

If $b$ is algebraic over $\FF_q$ then $K=\FF_q(b)$ is finite so
$[a,b)_{p^n}=0$ follows from $\Br (K)=\{0\}$. So we may assume that
$b$ is transcendental over $\FF_q$. It follows that $K$ is a global
field. Then we have the exact sequence
$$0\mapsto\Br (K)\to\bigoplus_{v\in\Omega_K}\Br (K_v)\to\QQ/\ZZ\to
0.$$
Here $\Omega_K$ is the set of all places of $K$ and the first map is
given by the localizations, $\xi\mapsto (\xi_v)_{v\in\Omega_K}$,
i.e. $[A]\mapsto ([A\otimes_KK_v])_{v\in\Omega_K}$ for any c.s.a.
$A$ over $K$. The second map is given by
$(\xi_v)_{v\in\Omega_K}\mapsto\sum_{v\in\Omega_K}inv_v(\xi_v)$, where
$inv_v:\Br (K_v)\xrightarrow\sim\QQ/\ZZ$ is the local invariant of the
Brauer group. Then for any $\xi\in\Br (K)$ we have $\xi =0$ iff
$\xi_v=0$ $\forall v\in\Omega_K$. If $v_0\in\Omega_K$ and $\xi_v=0$
only holds for $v\in\Omega_K\setminus\{ v_0\}$ then
$0=\sum_{v\in\Omega_K}inv_v(\xi_v)=inv_{v_0}(\xi_{v_0})$ so
$\xi_{v_0}=0$ as well. Therefore in order to prove that $\xi =0$ it
suffices to prove that $\xi_v=0$ holds for all but one value of
$v\in\Omega_K$.

If $v\in\Omega_K$ then we denote by ${\mathcal O}_v$ the ring of
integers from $K_v$, by ${\mathfrak p}_v$ the prime ideal and by
${\mathcal O}_v^\times ={\mathcal O}_v\setminus{\mathfrak p}_v$ the
group of units.

For every monic irreducible $f\in\FF_q[b]$ we have the place $v_f$ of
$K$ corresponding to the prime ideal $(f)$ of $\FF_q[b]$. Besides
these places, we have the place $v_\infty$ corresponding to the norm
$|\cdot |_\infty$, given by $|g/h|_\infty =q^{\deg g-\deg h}$ for every
$g,h\in\FF_q[b]\setminus\{0\}$.

Let $\xi =[a,b)_{p^n}\in\Br (K)$. Then for every $v\in\Omega_K$ we
have $\xi_v=[a,b)_{v,p^n}:=[a,b)_{K_v,p^n}$. So in order to prove that
$[a,b)_{p^n}=0$ it is enough to prove that $[a,b)_{v,p^n}=0$ holds for
every place $v$, except $v=v_\infty$. We have two cases.

If $v=v_f$ for some monic irreducible $f\in\FF_q[b]$, $f\neq b$, then
$b\in{\mathcal O}_v^\times$ and the entries of $a$ belong to
$b\FF_q[b]\sbq{\mathcal O}_v$, which, by [T, Corollary 2.1], implies
that $L_v:=K_v(\wp^{-1}(a))$ is an unramified extension of
$K_v$. Since $L_v/K_v$ is unramified and $b\in{\mathcal O}_v^\times$,
we have $b\in{\rm N}_{L_v/K_v}(L_v^\times )$ and so
$[a,b)_{v,p^n}=0$.

If $v=v_b$ then the entries of $a$ belong to $b\FF_q[b]\sbq{\mathfrak
p}_v$. By [T, Proposition 6.1], this implies that $a\in\wp
(W_n(K_v))$, so again $[a,b)_{v,p^n}=0$. 

\section{The symbols $((\cdot,\cdot ))_{p^m,p^n}$ and $((\cdot,\cdot
 ))_{p^m,p^n}$}

From now on we make the convention that $[0,0)_{p^n}=0$. 

\bdf For any field $K$ of characteristic $p$ and $n\geq 1$ we define
the symbol
$$((\cdot,\cdot ))_{p^n}:W_n(K)\times W_n(K)\to\Brr{p^n}(K)$$
as follows. If $a=(a_0,\ldots,a_{n-1}),\, b=(b_0,\ldots,b_{n-1})\in
W_n(K)$ then
$$((a,b))_{p^n}:=\sum_{j=0}^{n-1}[F^ja[b_j],b_j)_{p^n}.$$
By our convention, if $b_j=0$ then
$[F^ja[b_j],b_j)_{p^n}=[0,0)_{p^n}=0$ so the terms with $b_j=0$
should be ignored in the sum above.

In particular, if $a=0$ or $b=0$ then $((a,b))_{p^n}=0$.
\edf

\bff {\bf Remarks}

(1) If $a\in W_n(K)$, $b\in K$ then all but the first term in the
definition of $((a,[b]))_{p^n}$ are zero and we have
$((a,[b]))_{p^n}=[a[b],b)_{p^n}$.

Thus $[\cdot,\cdot )_{p^n}$ writes in terms of $((\cdot,\cdot
))_{p^n}$ as $[a,b)_{p^n}=((a[b]^{-1},[b]))_{p^n}$.

(2) If $n=1$ then $((\cdot,\cdot ))_p$ is defined by
$((a,b))_p=[ab,b)_p$.
\eff

\blm (i) With the notations from Definition 1, we have
$$((a,b))_{p^n}=\sum_{i,j=0}^{n-1}[[a_i^{p^j}b_j^{p^i}],b_j^{p^i})_{p^n}.$$

(ii) If $a,b\in K$ and $k,l\geq 0$ then
$((V^k[a],V^l[b]))_{p^n}=[[a^{p^l}b^{p^k}],b^{p^k})_{p^n}$.
\elm
\pf (i) We prove that $[F^ja[b_j],b_j)_{p^n}=
\sum_{i=0}^{n-1}[[a_i^{p^j}b_j^{p^i}],b_j^{p^i})_{p^n}$ for every
$j$. Then our result follows by summation over $j$. 

If $b_j=0$ this is just $0=0$ so we may assume that $b_j\neq 0$. We
have
$$F^ja[b_j]=(a_0^{p^j},\ldots,a_{n-1}^{p^j})[b_j]=
(a_0^{p^j}b_j,\ldots,a_{n-1}^{p^j}b_j^{p^{n-1}})=
\sum_{i=0}^{n-1}V^i[a_i^{p^j}b_j^{p^i}].$$
It follows that
$$[F^ja[b_j],b_j)_{p^n}=
\sum_{i=0}^{n-1}[V^i[a_i^{p^j}b_j^{p^i}],b_i)_{p^n}=
\sum_{i=0}^{n-1}[[a_i^{p^j}b_j^{p^i}],b_j^{p^i})_{p^n}.$$
(We have $[V^ia,b)_{p^n}=p^i[a,b)_{p^n}=[a,b^{p^i})_{p^n}$.) Hence the
conclusion.

(ii) Assume first that $k,l\leq n-1$. We use (i) in the case when
$a_i=0$ for $i\neq k$ and $b_j=0$ if $j\neq l$, i.e. when $a=V^k[a_k]$
and $b=V^l[b_l]$. It follows that
$[[a_i^{p^j}b_j^{p^i}],b_j^{p^i})_{p^n}=0$ for $(i,j)\neq (k,l)$ and
so $((V^k[a_k],V^l[b_l]))_{p^n}=((a,b))_{p^n}=
[[a_k^{p^l}b_l^{p^k}],b_l^{p^k})_{p^n}$. When we drop the indices $k$
and $l$ we get our result.

If $k\geq n$ or $l\geq n$ then $V^k[a]$ or $V^l[b]=0$ so
$((V^k[a],V^l[b]))_{p^n}=0$. If $k\geq n$ then $b^{p^k}$ is a
$p^n$-power so $[[a^{p^l}b^{p^k}],b^{p^k})_{p^n}=0$ and we are
done. Similarly, if $l\geq n$ then
$[[a^{p^l}b^{p^k}],a^{p^l})_{p^n}=0$. But by Lemma 2.2 we also have
$0=[[a^{p^l}b^{p^k}],a^{p^l}b^{p^k})_{p^n}=
[[a^{p^l}b^{p^k}],a^{p^l})_{p^n}+[[a^{p^l}b^{p^k}],b^{p^k})_{p^n}$ so
again $[[a^{p^l}b^{p^k}],b^{p^k})_{p^n}=0$. \qed

\bpr The symbol $((\cdot,\cdot))_{p^n}$ has the following properties.

(i) $((a,b))_{p^n}=((a+F^nc,b+F^nd))_{p^n}$ $\forall a,b,c,d\in
W_n(K)$.

(ii) $((\cdot,\cdot ))_{p^n}$ is bilinear.

(iii) $((a,b))_{p^n}=-((b,a))_{p^n}$ $\forall a,b\in W_n(K)$,
i.e. $((\cdot,\cdot))_{p^n}$ is skew-symmetric.

(iv) $((a,bc))_{p^n}+((b,ac))_{p^n}+((c,ab))_{p^n}=0$ $\forall
a,b,c\in W_n(K)$.
\epr
\pf (iii) If $a=(a_0,\ldots,a_{n-1})$ and $b=(b_0,\ldots,b_{n-1})$
then by Lemma 3.2(i) we have $((a,b))_{p^n}+((b,a))_{p^n}=
\sum_{i,j=0}^{n-1}[[a_i^{p^j}b_j^{p^i}],b_j^{p^i})_{p^n}+
\sum_{i,j=0}^{n-1}[[a_i^{p^j}b_j^{p^i}],a_i^{p^j})_{p^n}$. But by
Lemma 2.2 for every $i,j$ we have
$0=[[a_i^{p^j}b_j^{p^i}],a_i^{p^j}b_j^{p^i})_{p^n}=
[[a_i^{p^j}b_j^{p^i}],a_i^{p^j})_{p^n}+
[[a_i^{p^j}b_j^{p^i}],b_j^{p^i})_{p^n}$. Hence
$((a,b))_{p^n}+((b,a))_{p^n}=0$.

(ii) The linearity of $((\cdot,\cdot ))_{p^n}$ in the first variable
follows directly from the definition. Then the linearity in the second
variable will follow from the skew-symmetry, which we have already
proved.

(i) Since $((\cdot,\cdot ))_{p^n}$ is bilinear it suffices to prove
that $((F^na,b))_{p^n}=((a,F^nb))_{p^n}=0$ $\forall a,b\in
W_n(K)$. By using the formula $[a,b^{p^n})_{p^n}=0$ we get
$((a,F^nb))_{p^n}=\sum_{j=0}^{n-1}[F^ja[b_j^{p^n}],b_j^{p^n})_{p^n}=0$. Then
$((F^na,b))_{p^n}=0$ follows from the skew-symmetry.

(iv) We must prove that the map $f:W_n(K)^3\to\Brr{p^n}(K)$ given
by $(a,b,c)\mapsto ((a,bc))_{p^n}+((b,ac))_{p^n}+((c,ab))_{p^n}$ is
identically zero. Now $((\cdot,\cdot ))_{p^n}$ is bilinear so $f$ is
linear in each variable. Since $(W_n(K),+)$ is generated by $S=\{
V^i[a]\, :\, a\in\kk,\, 0\leq i\leq n-1\}$ it suffices to prove that
$f(a,b,c)=0$ when $a,b,c\in S$. So we must prove that
$f(V^i[a],V^j[b],V^k[c])=0$ $\forall a,b,c\in\kk$, $0\leq i,j,k\leq
n-1$. We have $V^j[b]V^k[c]= V^{j+k}(F^k[b]F^j[c])=
V^{j+k}[b^{p^k}c^{p^j}]$. By Lemma 3.2(ii)
$((V^i[a],V^j[b]V^k[c]))_{p^n}= -((V^j[b]V^k[c],V^i[a]))_{p^n}$ writes
as
\begin{multline*}
-((V^{j+k}[b^{p^k}c^{p^j}],V^i[a]))_{p^n}=
-[[(b^{p^k}c^{p^j})^{p^i}a^{p^{j+k}}],a^{p^{j+k}})_{p^n}\\
=-[[a^{p^{j+k}}b^{p^{i+k}}c^{p^{i+j}}],a^{p^{j+k}})_{p^n}.
\end{multline*}
Similarly for the remaining two terms of
$f(V^i[a],V^j[b],V^k[c])$. Hence
\begin{multline*}
f(V^i[a],V^j[b],V^k[c])=
-[[a^{p^{j+k}}b^{p^{i+k}}c^{p^{i+j}}],a^{p^{j+k}})_{p^n}
-[[a^{p^{j+k}}b^{p^{i+k}}c^{p^{i+j}}],b^{p^{i+k}})_{p^n}\\
-[[a^{p^{j+k}}b^{p^{i+k}}c^{p^{i+j}}],c^{p^{i+j}})_{p^n}=
-[[a^{p^{j+k}}b^{p^{i+k}}c^{p^{i+j}}],a^{p^{j+k}}b^{p^{i+k}}c^{p^{i+j}})_{p^n},
\end{multline*}
which is zero by Lemma 2.2. \qed

Properties (i)-(iii) summarize as follows.

\bco $((\cdot,\cdot ))_{p^n}$ is a bilinear skew-symmetric map defined
as
$$((\cdot,\cdot ))_{p^n}:W_n(K)/F^n(W_n(K))\times
W_n(K)/F^n(W_n(K))\to\Brr{p^n}(K).$$
Note that $F^n(W_n(K))$ also writes as $W_n(K^{p^n})$. 
\eco

\blm For every ring $R$ there is a group isomorphism
$$(R\otimes R)/\langle a\otimes bc-ab\otimes c-ac\otimes b\, :\,
a,b,c\in R\rangle\to\Omega^1(R)$$
given by $x\otimes y\mapsto x\diff y$.
\elm
\pf $\Omega^1(R)$ is the $R$-module generated by $\diff a$ with $a\in
R$, subject to $\diff (a+b)=\diff a+\diff b$ and $\diff (ab)=a\diff
b+b\diff a$ $\forall a,b\in R$. Then $\Omega^1(R)$ writes as $M/N$,
where $M$ is the $R$-module generated by $\diff a$ with $a\in R$
subject to $\diff (a+b)=\diff a+\diff b$ and $N$ is the $R$-submodule
of $M$ generated by $\diff (ab)-a\diff b-b\diff a$, with $a,b\in R$.

We claim that there is a group isomorphism $f:R\otimes R\to M$ given
by $x\otimes y\mapsto x\diff y$. The existence of $f$ defined this way
follows from the fact that the map $R\times R\to M$ given by $(x,y)\to
x\diff y$ is bilinear. (In $M$ we have $(a+b)\diff c=a\diff c+b\diff
c$ and $a\diff (b+c)=a(\diff b+\diff c)=a\diff b+a\diff c$.)
Conversely, we regard $R\otimes R$ as an $R$-module by defining
$x\alpha :=(x\otimes 1)\alpha$ $\forall x\in R,\, \alpha\in R\otimes
R$ and we define a morphism of $R$-modules $g:M\to R\otimes R$ by
$\diff x\mapsto 1\otimes x$. This is well defined because the
relations among generators in $M$, $\diff (a+b)=\diff a+\diff b$, are
preserved by $g$. (We have $1\otimes (a+b)=1\otimes a+1\otimes b$.)
Now $g(x\diff y)=xg(\diff y)=x(1\otimes y)=(x\otimes 1)(1\otimes
y)=x\otimes y$. It follows that $f$ and $g$ are inverse to each other
group isomorphisms.

Then $f$ induces a group isomorphism $(R\otimes R)/g(N)\to
M/N=\Omega^1(R)$, given by $x\otimes y\to x\diff y$. Now, as an
$R$-module, $N$ is generated by $\diff (bc)-b\diff c-c\diff b$, with
$b,c\in R$. As a group, it will be generated by $a(\diff (bc)-b\diff
c-c\diff b)=a\diff (bc)-ab\diff c-ac\diff b$, with $a,b,c\in R$. It
follows that $g(N)$ is the group generated by $g(a\diff (bc)-ab\diff
c-ac\diff b)=a\otimes bc-ab\otimes c-ac\otimes b$, with $a,b,c\in
R$. Hence the conclusion. \qed

\bpr There is a group morphism
$\alpha_{p^n}:\Omega^1(W_n(K))/\diff W_n(K)\to\Brr{p^n}(K)$ given by
$a\diff b\mapsto ((a,b))_{p^n}$.

In particular, if $n=1$ then $((a,b))_p=[ab,b)_p$ (see Remark 3.1(2))
so we recover the original definition of $\alpha_p$ from the
introduction.
\epr
\pf For convenience, we write $((\cdot,\cdot ))$ instead of
$((\cdot,\cdot ))_{p^n}$. By Proposition 3.3(ii) $((\cdot,\cdot
)):W_n(K)\times W_n(K)\to\Brr{p^n}(K)$ is bilinear so there is
a group morphism $f:W_n(K)\otimes W_n(K)\to\Brr{p^n}(K)$ given
by $a\otimes b\mapsto ((a,b))$. By Proposition 3.3(iii) and (iv) for
every $a,b,c\in W_n(K)$ we have $f(a\otimes bc-ab\otimes c-ac\otimes
b)=((a,bc))-((ab,c))-((ac,b))=((a,bc))+((c,ab))+((b,ac))=0$. So $f$
can be defined on
$$(W_n(K)\otimes W_n(K))/\langle a\otimes bc-ab\otimes
c-ac\otimes b\, :\, a,b,c\in W_n(K)\rangle,$$
which, by Lemma 3.5, is isomorphic to $\Omega^1(W_n(K))$, via
$a\otimes b\mapsto a\diff b$. Then we get a group morphism
$\alpha_{p^n}:\Omega^1(W_n(K))\to\Brr{p^n}(K)$ given by $a\diff
b\mapsto f(a\otimes b)=((a,b))$. But for every $a,b\in W_n(K)$ we
have $\alpha_{p^n}(\diff (ab))=\alpha_{p^n}(a\diff b+b\diff
a)=((a,b))+((b,a))=0$. In particular, if $b=1$ we get
$\alpha_{p^n}(\diff a)=0$ $\forall a\in W_n(K)$. Hence
$\alpha_{p^n}$ is defined in fact on $\Omega^1(W_n(K))/\diff
W_n(K)$. \qed

{\bf Remark} Proposition 3.6 is a consequence of Proposition
3.3(ii)-(iv). But in fact we have equivalence. Indeed, (ii) follows
from $\alpha_{p^n}(a\diff (b+c))=\alpha_{p^n}(a\diff
b)+\alpha_{p^n}(a\diff c)$ and $\alpha_{p^n}((a+b)\diff
c)=\alpha_{p^n}(a\diff c)+\alpha_{p^n}(b\diff c)$. For (iii) we have
$\alpha_{p^n}(a\diff b+b\diff a)=\alpha_{p^n}(\diff (ab))=0$,
i.e. $((a,b))+((b,a))=0$. And for (iv) we have $\alpha_{p^n}(a\diff
(bc))=\alpha_{p^n}(ab\diff c+ac\diff b)$,
i.e. $((a,bc))=((ab,c))+((ac,b))$. Together with $((ab,c))=-((c,ab))$
and $((ac,b))=-((b,ac))$, this implies
$((a,bc))+((b,ac))+((c,ab))=0$.

\bpr The Frobenius and Verschiebung maps are adjoint:
$$((Fa,b))_{p^n}=((a,Vb))_{p^n}\quad\text{and}\quad
((Va,b))_{p^n}=((a,Fb))_{p^n}\quad\forall a,b\in W_n(K).$$
\epr
\pf Since both maps $(a,b)\mapsto ((Fa,b))_{p^n}$ and
$(a,b)\mapsto ((a,Vb))_{p^n}$ are bilinear and $(W_n(K),+)$ is
generated by $S=\{ V^i[a]\, :\, a\in\kk,\, 0\leq i\leq n-1\}$, it
suffices to prove that $((Fa,b))_{p^n}=((a,Vb))_{p^n}$ for $a,b\in
S$. So we must prove that
$((FV^i[a],V^j[b]))_{p^n}=((V^i[a],VV^j[b]))_{p^n}$ $\forall a,b\in\kk$,
$0\leq i,j\leq n-1$. We use Lemma 3.2(ii) and we get
\begin{multline*}
((FV^i[a],V^j[b]))_{p^n}= ((V^i[a^p],V^j[b]))_{p^n}=
[[(a^p)^{p^j}b^{p^i}],b^{p^i})_{p^n}\\
=[[a^{p^{j+1}}b^{p^i}],b^{p^i})_{p^n}=((V^i[a],V^{j+1}[b]))_{p^n}.
\end{multline*}

The second statement follows from the first by the skew-symmetry. \qed

\bpr If $a,b\in W_n(K)$ then
$((a,b))_{p^n}=((Va,Vb))_{p^{n+1}}$.

More generally, if $m\geq n$ then
$((a,b))_{p^n}=((V^{m-n}a,V^{m-n}b))_{p^m}$.
\epr
\pf Let $b=(b_0,\ldots,b_{n-1})$. Then $Vb=(0,b_0,\ldots,b_{n-1})$ so,
by definition,
$$((Va,Vb))_{p^{n+1}}= [0,0)_{p^{n+1}}+
\sum_{j=0}^{n-1}[F^{j+1}Va[b_j],b_j)_{p^{n+1}}.$$
But $F^{j+1}Va[b_j]=V(F^{j+1}a)[b_j]=V(F^{j+1}aF
[b_j])$. It follows that

$$[F^{j+1}Va[b_j],b_j)_{p^{n+1}}= [V(F^{j+1}aF[b_j]),b_j)_{p^{n+1}}=
[F^{j+1}aF[b_j],b_j)_{p^n}= [F^ja[b_j],b_j)_{p^n}.$$

(We used the formulas $[Va,b)_{p^{n+1}}=[a,b)_{p^n}$ and
$[Fa,b)_{p^n}=[a,b)_{p^n}$.) 

Hence
$((Va,Vb))_{p^{n+1}}=\sum_{j=0}^{n-1}[F^ja[b_j],b_j)_{p^n}=((a,b))_{p^n}$.
\qed

\bco If $m\geq n$ then for any $a,b\in W(K)$ we have
$((a,b))_{p^n}=p^{m-n}((a,b))_{p^m}$. Explicitly, 
$$(((a_0,\ldots,a_{n-1}),(b_0,\ldots,b_{n-1})))_{p^n}=
p^{m-n}(((a_0,\ldots,a_{m-1}),(b_0,\ldots,b_{m-1})))_{p^m}.$$
\eco
\pf By Propositions 3.8 and 3.7, we have
$((a,b))_{p^n}=((V^{m-n}a,V^{m-n}b))_{p^m}=((F^{m-n}V^{m-n}a,b))_{p^m}=
((p^{m-n}a,b))_{p^m}=p^{m-n}((a,b))_{p^m}$. \qed

\bpr $(\cdot,\cdot )_{p^n}$ is antisymmetric. 
\epr
\pf Since $(\cdot,\cdot )_{p^n}$ is skew-symmetric we have
$2((a,a))_{p^n}=((a,a))_{p^n}+((a,a))_{p^n}=0$. If $p>2$ then also
$p^n((a,a))_{p^n}=0$. Since $(2,p^n)=1$ we get $((a,a))_{p^n}=0$. If
$p=2$ then by Corollary 3.9 we have
$((a,a))_{2^n}=2((a,a))_{2^{n+1}}=0$. \qed

\bdf For $m,n\geq 1$ we define the symbol
$$((\cdot,\cdot ))_{p^m,p^n}:W_m(K)\times W_n(K)\to\Br (K)$$
by $((a,b))_{p^m,p^n}=((V^{l-m}a,V^{l-n}b))_{p^l}$ for any $l\geq m,n$.

In particular, if $m=n$ we may take $l=n$ and we have $((\cdot,\cdot
))_{p^n,p^n}=((\cdot,\cdot ))_{p^n}$.
\edf

\bpr (i) $((\cdot,\cdot))_{p^m,p^n}$ is well defined.

(ii) $((a,b))_{p^m,p^n}=((a+F^nc,b+F^md))_{p^m,p^n}$ $\forall
a,c\in W_m(K)$, $b,d\in W_n(K)$.

(iii) $((\cdot,\cdot))_{p^m,p^n}$ is bilinear.

(iv) $((a,b))_{p^m,p^n}=-((b,a))_{p^n,p^m}$ $\forall a\in W_m(K)$,
$b\in W_n(K)$.
\epr
\pf (i) We must prove that the formula for $((a,b))_{p^m,p^n}$ from
Definition 2 is independent of the choice of $l$. Assume that
$l'\geq l\geq m,n$. Then by Proposition 3.8 we have
$((V^{l-m}a,V^{l-n}b))_{p^l}=
((V^{l'-l}V^{l-m}a,V^{l'-l}V^{l-n}b))_{p^{l'}}= ((V^{l'-m}a,V^{l'-n}b))_{p^{l'}}$.

(iii) follows from the bilinearity of $((\cdot,\cdot))_{p^l}$ and the
fact that the maps $a\mapsto V^{l-m}a$ and $b\mapsto V^{l-n}b$ are
linear.

(ii) Since $((\cdot,\cdot))_{p^m,p^n}$ is bilinear it suffices to
prove that $((a,F^mb))_{p^m,p^n}=((F^na,b))_{p^m,p^n}=0$
$\forall a\in W_m(K)$, $b\in W_n(K)$. If $l\geq m,n$ then
\begin{multline*}
((a,F^mb))_{p^m,p^n}= ((V^{l-m}a,V^{l-n}F^mb))_{p^l}\\
=((V^mV^{l-m}a,V^{l-n}b))_{p^l}= ((0,V^{l-n}b))_{p^l}=0.
\end{multline*}
(Here we used the adjoint property of $F$ and $V$ and the fact that
$V^l\equiv 0$ on $W_l(K)$.) The proof of $((F^na,b))_{p^m,p^n}=0$ is
similar.

(iv) Follows from the definition of $((\cdot,\cdot ))_{p^m,p^n}$ and
the skew-symmetry of $((\cdot,\cdot ))_{p^l}$. \qed

\bco $((\cdot,\cdot ))_{p^m,p^n}$ is a bilinear map defined as
$$((\cdot,\cdot ))_{p^m,p^n}:W_m(K)/F^n(W_m(K))\times
W_n(K)/F^m(W_n(K))\to\Brr{p^k}(K),$$
where $k=\min\{m,n\}$.

Note that $F^n(W_m(K))$ and $F^m(W_n(K))$ also write as $W_m(K^{p^n})$
and $W_n(K^{p^m})$.
\eco
\pf By Proposition 3.11(ii), $((a,b))_{p^m,p^n}$ depends only on
$a\mod F^n(W_m(K))$ and $b\mod F^m(W_n(K))$. This justifies the
new domain for $((\cdot,\cdot ))_{p^m,p^n}$. The fact that the image
of  $((\cdot,\cdot ))_{p^m,p^n}$ is in $\Brr{p^k}(K)$ follows from
bilinearity of $((\cdot,\cdot ))_{p^m,p^n}$ and the fact that
$W_m(K)/F^n(W_m(K))$ and $W_n(K)/F^m(W_n(K))$
are $p^k$-torsion. Indeed, $FV=VF=p$ so $p^mW_m(K)\sbq
V^m(W_m(K))=\{0\}$ and $p^nW_m(K)\sbq F^n(W_m(K))$. Thus
$W_m(K)/F^n(W_m(K))$ is killed by both $p^m$ and $p^n$ and
so by $p^k$. Similarly for $W_n(K)/F^m(W_n(K))$. \qed

\bpr The Frobenius and Verschiebung maps are adjoint:

$((Fa,b))_{p^m,p^n}=((a,Vb))_{p^m,p^n}=((a,b))_{p^m,p^{n-1}}$ 

$((Va,b))_{p^m,p^n}=((a,Fb))_{p^m,p^n}=((a,b))_{p^{m-1},p^n}$

Here we make the convention that $((a,b))_{p^m,p^n}=0$ if $m$ or
$n=0$.

More generally, $\forall a,b\in W(K)$ we have
$$((F^iV^ja,F^kV^lb))_{p^m,p^n}=\begin{cases}
((a,b))_{p^{m-j-k},p^{n-i-l}}&\text{if }m>j+k,\, n>i+l\\
0&\text{otherwise}
\end{cases}.$$
\epr
\pf Let $N\geq m,n$. Since $F$ and $V$ are adjoint with respect to
$((\cdot,\cdot ))_{p^N}$ we have 
\begin{multline*}
((F^iV^ja,F^kV^lb))_{p^m,p^n}=
((V^{N-m}F^iV^ja,V^{N-n}F^kV^lb))_{p^N}\\
=((F^iV^{N-m+j}a,F^kV^{N-n+l}b))_{p^N}=
((V^kV^{N-m+j}a,V^iV^{N-n+l}b))_{p^N}.
\end{multline*}
If $m>j+k$, $n>i+l$ then $N\geq m-j-k,n-i-l\geq 1$ so, by definition,
$((V^{N-m+j+k}a,V^{N-n+l+i}b))_{p^N}=((a,b))_{p^{m-j-k},p^{n-i-l}}$. If
$m\leq j+k$ or $n\leq i+l$ then $N-m+j+k$ or $N-n+l+i\geq N$ so
$V^{N-m+j+k}a$ or $V^{N-n+l+i}b=0$ in $W_N(K)$. (On $W_N(K)$ we have
$V^N\equiv 0$.) Hence $((V^{N-m+j+k}a,V^{N-n+l+i}b))_{p^N}=0$. \qed

\bff{\bf Remark.} In short notation we may write
$W_m(K)/F^n(W_m(K))$ as $W_m(K)/(F^n)$, where by
$(F^n)$ we mean the image of $F^n$ on $W_m(K)$. In the same
short notation $W_m(K)=W(K)/V^m(W(K))$ may be written as
$W(K)/(V^m)$. Hence in the short notation we have
$W_m(K)/F^n(W_m(K))=W(K)/(V^m,F^n)$, where
$(V^m,F^n)$ is the subgroup of $W(K)$ generated by the images of
$V^m$ and $F^n$. Similarly,
$W_n(K)/F^m(W_n(K))=W(K)/(V^n,F^m)$ so the domain of
$((\cdot,\cdot ))_{p^m,p^n}$ can be written as
$W(K)/(V^m,F^n)\times W(K)/(V^n,F^m)$. 

One can see that $F$ and $V$ switch roles in $W(K)/(V^m,F^n)$ and
$W(K)/(V^n,F^m)$. This is explained by the fact that $F$ and $V$
are adjoint with respect to $((\cdot,\cdot ))_{p^m,p^n}$ so if
$P\in\ZZ [X,Y]$ then
$((P(V,F)a,b))_{p^m,p^n}=((a,P(F,V)b))_{p^m,p^n}$ $\forall
a,b\in W(K)$. In particular, $((P(V,F)a,b))_{p^m,p^n}=0$ $\forall
a,b\in W(K)$ iff $((a,P(F,V)b))_{p^m,p^n}=0$ $\forall a,b\in
W(K)$. Two polynomials $P$ with this property are $X^m$ and $Y^n$.
\eff

\bpr (i) If $a=(a_0,\ldots,a_{m-1})\in W_m(K)$,
$b=(b_0,\ldots,b_{n-1})\in W_n(K)$ then
$$((a,b))_{p^m,p^n}=\sum_{j=0}^{n-1}[F^{m+j}aF^n[b_j],b_j)_{p^m}.$$
(Same as in Definition 1, the terms with $b_j=0$ should be ignored.)

(ii) $[F^{m+j}aF^n[b_j],b_j)_{p^m}= ((a,V^j[b_j]))_{p^m,p^n}=
((a,[b]_j))_{p^m,p^{n-j}}$ $\forall 0\leq j\leq n-1$.
\epr
\pf (i) In Definition 2 we take $l=m+n$ and we get
$((a,b))_{p^m,p^n}=((V^na,V^mb))_{p^{m+n}}$. For $0\leq j\leq n-1$ the
$m+j$ entry of $V^mb$ is $b_j$ and all the other entries are $0$. By
Definition 1 we have 
$$((V^na,V^mb))_{p^{m+n}}=
\sum_{j=0}^{n-1}[F^{m+j}V^na[b_j],b_j)_{p^{m+n}}.$$
But $F^{m+j}V^na[b_j]= V^n(F^{m+j}a)[b_j]=
V^n(F^{m+j}aF^n[b_j])$. It follows that 
$$[F^{m+j}V^na[b_j],b_j)_{p^{m+n}}=
[V^n(F^{m+j}aF^n[b_j]),b_j)_{p^{m+n}}=
[F^{m+j}aF^n[b_j],b_j)_{p^m}.$$
Hence the conclusion.

(ii) The $j$ entry of $V^j[b_j]$ is $b_j$ and all the other entries
are $0$. Then by (i) we have
$((a,V^j[b_j]))_{p^m,p^n}=[F^{m+j}aF^n[b_j],b_j)_{p^m}$.

The relation $((a,V^j[b_j]))_{p^m,p^n}= ((a,[b_j]))_{p^m,p^{n-j}}$
follows from Proposition 3.13. \qed

\bff{\bf Remarks}

(1) Since $[Fa,b)_{p^m}=[a,b)_{p^m}$, we can simplify $F^k$,
where $k=\min\{ m+j,n\}$, in the formula
$[F^{m+j}aF^n[b_j],b_j)_{p^m}$. So
$[F^{m+j}aF^n[b_j],b_j)_{p^m}=[F^{m-n+j}a[b_j],b_j)_{p^m}$ if
$m+j\geq n$ and $=[aF^{n-m-j}[b_j],b_j)_{p^m}$ if $m+j<n$.

So if $m\geq n$ then $((a,b))_{p^m,p^n}=
\sum_{j=0}^{n-1}[F^{m-n+j}a[b_j],b_j)_{p^m}$, while if $m<n$ then
$((a,b))_{p^m,p^n}= \sum_{j=0}^{n-m-1}[aF^{n-m-j}[b_j],b_j)_{p^m}+
\sum_{j=n-m}^{n-1}[F^{m-n+j}a[b_j],b_j)_{p^m}$.

In particular, when $m=n$ we recover Definition 1 for $((\cdot,\cdot
))_{p^n}=((\cdot,\cdot ))_{p^n,p^n}$.

(2) A priori, $[F^{m+j}aF^n[b_j],b_j)_{p^m}$ is
$p^m$-torsion. But in fact, since it writes as
$((a,[b_j]))_{p^m,p^{n-j}}$, it is $p^{\min\{ m,n-j\}}$-torsion. In
particular, when $m=n$ the term $[F^ja[b_j],b_j)_{p^n}$ from
Definition 1 writes as $((a,[b_j]))_{p^n,p^{n-j}}$ so it is
$p^{n-j}$-torsion.

(3) We have $((a,b))_{p^m,p^n}=-((b,a))_{p^n,p^m}$ so 
$$((a,b))_{p^m,p^n}=
\sum_{j=0}^{n-1}[F^{m+j}aF^n[b_j],b_j)_{p^m}=
-\sum_{i=0}^{m-1}[F^{n+i}bF^m[a_i],a_i)_{p^n}.$$
In terms of c.s.a., $((a,b))_{p^m,p^n}=[A]=[B]$, where
$A=\bigotimes_{j=0}^{n-1}A_{[F^{m+j}aF^n[b_j],b_j)_{p^m}}$ and
$B=\bigotimes_{i=0}^{m-1}A_{[F^{n+i}bF^m[a_i],a_i)_{p^n}}^{op}$. Thus
$A$ writes as the tensor product of $n$ c.s.a. of degree $p^m$ and $B$
as the tensor product of $m$ c.s.a. of degree $p^n$. Hence $\deg
A=\deg B=p^{mn}$. Therefore the Schur index of $((a,b))_{p^m,p^n}$ is
at most $p^{mn}$. Philippe Gille raised the question whether this
upper bound can be achieved. The answer is yes. If
$F=\FF_p(a_0,\ldots,a_{m-1},b_0,\ldots,b_{n-1})$, where $a_i,b_j$ are
variables, and $a=(a_0,\ldots,a_{m-1})$, $b=(b_0,\ldots,b_{n-1})$ then
the Schur index of $((a,b))_{p^m,p^n}$ is $p^{mn}$. To prove this we
need a new way to describe $((\cdot,\cdot ))_{p^m,p^n}$ in terms of
c.s.a. given by generators and relations.
\eff

As seen in the introduction, the symbols $[\cdot,\cdot
)_{p^n}:W_n(K)\times\kk\to\Brr{p^n}(K)$ have a direct limit
$[\cdot,\cdot )_{p^\infty}:CW(K)\times\kk\to\Brr{p^\infty}(K)$. For
every $n\geq 1$ the canonical morphism $\psi_n:W_n(K)\to CW(K)$ is
given by $(a_0,\ldots,a_{n-1})\mapsto
(\ldots,0,0,a_0,\ldots,a_{n-1})$.

We can do the same with the symbols $((a,b))_{p^m,p^n}$ indexed over
$(\NN^*\times\NN^*,\leq )$, where $(m,n)\leq (m',n')$ if $m\leq m'$,
$n\leq n'$.

If $m\leq m'$, $n\leq n'$ then, by Proposition 3.13, for any $a\in
W_m(K)$, $b\in W_n(K)$ we have
$((a,b))_{p^m,p^n}=((V^{m'-m}a,V^{n'-n}b))_{p^{m'},p^{n'}}$. So we
have the commuting diagram
$$\begin{array}{ccc}W_m(K)\times W_n(K)&\xrightarrow{((\cdot,\cdot
))_{p^m,p^n}}&\Brr{p^k}(K)\\
\hskip -60pt V^{m'-m}\times V^{n'-n}\downarrow&{}&\downarrow\\
W_{m'}(K)\times W_{n'}(K)&\xrightarrow{((\cdot,\cdot
))_{p^{m'},p^{n'}}}&\Brr{p^{k'}}(K)
\end{array},$$
where $k=\min\{ m,n\}$, $k'=\min\{ m',n'\}$. So we have a map between
two directed systems. By taking direct limits we get a symbol
$((\cdot,\cdot ))_{p^\infty}:CW(K)\times CW(K)\to\Brr{p^\infty}(K)$. If
$a=(\ldots,a_{-1},a_0)$, $b=(\ldots,b_{-1},b_0)\in CW(K)$ with $a_i=0$
for $i\leq -m$ and $b_j=0$ for $j\leq -n$ then
$a=\psi_m((a_{-m+1}\ldots,a_0))$ and
$b=\psi_n((b_{-n+1}\ldots,b_0))$. So $((a,b))_{p^\infty}=
(((a_{-m+1}\ldots,a_0),(b_{-n+1}\ldots,b_0)))_{p^m,p^n}$.

Now $\{ (n,n)\, :\, n\in\NN^*\}$ is cofinal in $(\NN^*\times\NN^*,\leq
)$ so $((\cdot,\cdot ))_{p^\infty}$ can be regarded also as the direct
limit of $((\cdot,\cdot ))_{p^n,p^n}=((\cdot,\cdot ))_{p^n}$
only. Since $((\cdot,\cdot ))_{p^n}$ are bilinear and antisymmetric, so
is $((\cdot,\cdot ))_{p^\infty}$. 

Note that if $a\in CW(K^{p^\infty})$ then $a_i\in K^{p^\infty}\sbq
K^{p^n}$ $\forall i$. Hence $(a_{-m+1}\ldots,a_0)\in
W_m(K^{\p^n})=F^n(W_m(K))$. It follows that
$(((a_{-m+1}\ldots,a_0),(b_{-n+1}\ldots,b_0)))_{p^m,p^n}=0$,
i.e. $((a,b))_{p^\infty}=0$. Similarly, $((a,b))_{p^\infty}=0$ if
$b\in CW(K^{p^\infty})$. Since $((\cdot,\cdot ))_{p^\infty}$ is
bilinear, this implies that $((a,b))_{p^\infty}$ depends only on $a$
and $b\mod CW(K^{p^\infty})$. So the symbol can be defined as
$((\cdot,\cdot ))_{p^\infty}:CW(K)/CW(K^{p^\infty})\times
CW(K)/CW(K^{p^\infty})\to\Brr{p^\infty}(K)$. We get:

\bpr There is a bilinear antisymmetric symbol
$$((\cdot,\cdot ))_{p^\infty}:CW(K)/CW(K^{p^\infty})\times
CW(K)/CW(K^{p^\infty})\to\Brr{p^\infty}(K)$$
given for any $a=(\ldots,a_{-1},a_0)$, $b=(\ldots,b_{-1},b_0)\in
CW(K)$ with $a_i=0$ if $i\leq -m$ and $b_j=0$ if $j\leq -n$ by
$((a,b))_{p^\infty}=
(((a_{-m+1}\ldots,a_0),(b_{-n+1}\ldots,b_0)))_{p^m,p^n}$. 
\epr

{\bf Remark.} The symbols $[\cdot,\cdot )_{p^n}$ too satisfy the
adjoint property between Frobenius and Verschiebung. We have
$[Va,b)_{p^n}=[a,b)_{p^{n-1}}$, but also
$[a,Fb)_{p^n}=[a,b^p)_{p^n}=p[a,b)_{p^n}=[a,b)_{p^{n-1}}$. So
$[Va,b)_{p^n}=[a,Fb)_{p^n}=[a,b)_{p^{n-1}}$.

For the other adjoint property we need to define a Verschiebung map on
$\kk$. On $(W_n(F),+)$ the Verschiebung map is defined by the property
$FV=VF=p:=(x\mapsto px)$. On $(\kk,\cdot )$ we need a multiplicative
Verschiebung map $V^\times$, which should satisfy $FV^\times =VF^\times
=p:=(x\mapsto x^p)$. Obvious such map is the identity, $V^\times =1$.
Then we have $[Fa,b)_{p^n}=[a,V^\times b)_{p^n}=[a,b)_{p^n}$.


As we will see in a future paper, this is a particular case of a more
general result.

\section{The representation theorem}

In Proposition 3.6 we introduced the linear map
$\alpha_{p^n}:\Omega^1(W_n(K))\to\Brr{p^n}(K)$ given by $a\diff
b\mapsto ((a,b))_{p^n}$ . We are now able to prove that $\alpha_{p^n}$
is surjective and to find its kernel, thus to generalize the
result of [GS, Theorem 9.2.4] for $n=1$, which we mentioned in the
introduction. 

Note that in fact we already have the surjectivity. Indeed,
$\Ima\alpha_{p^n}$ is the subgroup of $\Brr{p^n}(K)$ generated by the
image of $((\cdot,\cdot ))_{p^n}$, which, by Remark 3.1(1), coincides
with the subgroup generated by the image of $[\cdot,\cdot
)_{p^n}$, which, by Proposition 1.2, is the whole $\Brr{p^n}(K)$.

\blm For every $a\in W_n(K)$, $b\in\kk$ we have
$[a,b)_{p^n}=\alpha_{p^n}(a\dlog [b])$.
\elm
\pf By Remark 3.1 we have
$[a,b)_{p^n}=((a[b]^{-1},[b]))_{p^n}=\alpha_{p^n}(a[b]^{-1}\diff
[b])=\alpha_{p^n}(a\dlog [b])$. \qed

\blm The following elements of $\Omega^1(W_n(K))$ belong to
$\ker\alpha_{p^n}$.

\begin{center}
$\diff a$, $Fa\diff b-a\diff Vb$, $Va\diff b-a\diff Fb$,
$a,b\in W_n(K)$.

$\wp (a)\dlog [b]=(Fa-a)\dlog [b]$, $a\in W_n(K)$, $b\in\kk$. 
\end{center}

In particular, $\wp ([a])\dlog [b]=([a^p]-[a])\dlog
[b]\in\ker\alpha_{p^n}$, $\forall a\in K$, $b\in\kk$.
\elm
\pf We have $\diff a\in\ker\alpha_{p^n}$ by Proposition 3.6. 

By Proposition 3.7 we have $\alpha_{p^n}(Fa\diff b-a\diff
Vb)=((Fa,b))_{p^n}-((a,Vb))_{p^n}=0$. Similarly for $Va\diff
b-a\diff Fb$.

By Lemma 4.1 we have $\alpha_{p^n}(\wp (a)\dlog
[b])=[\wp(a),b)_{p^n}=0$. \qed

\bdf We define $G_n:=\Omega^1(W_n(K))/M_n$, where 
$$M_n=\langle Fa\diff b-a\diff Vb\, :\, a,b\in W_n(K),\,\wp
([a])\dlog [b]\, :\, a\in W_n(K),\, b\in\kk\rangle.$$

By Lemma 4.2 we have $M_n\sbq\ker\alpha_{p^n}$. Therefore we can
regard $\alpha_{p^n}$ as being defined
$\alpha_{p^n}:G_n\to\Brr{p^n}(K)$.

We also define $G'_n:=\Omega^1(W_n(K))/M'_n$, where $M'_n\sbq M_n$,
$$M'_n=\langle Fa\diff b-a\diff Vb\, :\, a,b\in W_n(K)\rangle.$$
\edf	

We will prove by induction on $n$ that
$\alpha_{p^n}:\Omega^1(W_n(K))\to\Brr{p^n}(K)$ is surjective and
$M_n$ is its kernel so $\alpha_{p^n}:G_n\to\Brr{p^n}(K)$ is an
isomorphism.

\blm We have $M_1=\langle\diff a,\, \wp (a)\dlog b\, :\, a\in K,\,
b\in\kk\rangle$. 

Consequently, $\alpha_p:G_1\to\Brr p(K)$ is an isomorphism by [GS,
Theorem 9.2.4]. (See also \S 1.)
\elm
\pf In $W_1(K)=K$ we have $V\equiv 0$ and $[a]=a$ $\forall a\in K$. It
follows that $Fa\diff b-a\diff Vb=a^p\diff b-a\diff 0=\diff
a^pb$. Thus $\{Fa\diff b-a\diff Vb\, :\, a,b\in K\} =\{\diff a\, :\,
a\in K\}$. Also $\wp ([a])\dlog [b]=\wp (a)\dlog b$ $\forall a\in K$,
$b\in\kk$. Hence the conclusion. \qed

Note that we used only a minimal set of generators for $M_n$. In fact
$M_n$ contains all elements of $\ker\alpha_{p^n}$ we know so far.

\blm All elements of $\ker\alpha_{p^n}$ from Lemma 4.2 also belong to
$M_n$.

With the exception of $\wp (a)\dlog [b]$, they also belong to $M'_n$.
\elm
\pf In $G'_n$ we have $Fa\diff b=a\diff Vb$. More generally,
$F^ka\diff b=a\diff V^kb$. So $\diff a=F^n1\diff a=1\diff
V^na=0$. (We have $V^n\equiv 0$ in $W_n(K)$ so $V^na=0$.)
Consequently, $a\diff b+b\diff a=\diff (ab)=0$ so $a\diff b=-b\diff
a$. Hence $Va\diff b=-b\diff Va=-Fb\diff a=a\diff Fb$. The relations
$\diff a=0$ and $Va\diff b=a\diff Fb$, which hold in $G'_n$, imply
$\diff a,Va\diff b-a\diff Fb\in K'_n\sbq K_n$. Note that $\diff a=0$
and $Va\diff b=a\diff Fb$ also hold in $G_n$.

We are left to prove that $\wp a\dlog [b]\in M_n$, i.e. that $\wp
a\dlog [b]=0$ in $G_n$. The map $f:W_n(K)\to
G_n$, given by $a\mapsto\wp a\dlog [b]$, is linear. We must prove that
$f\equiv 0$. The group $(W_n(K),+)$ is generated by $[a]$, with
$a\in K$, and $Va$, where $a$ is a Witt vector, so it suffices to
prove that $f$ vanishes at these generators. We have $\wp ([a])\dlog
[b]\in K_n$ so $f([a])=0$. In $G_n$ we have $Va\dlog
[b]=Va[b]^{-1}\diff [b]=V(aF[b]^{-1})\diff [b]=aF[b]^{-1}\diff
F[b]=a\dlog F[b]$. But $\dlog F[b]=\dlog [b]^p=p\dlog [b]$. It follows
that  $Va\dlog [b]=pa\dlog [b]=FVa\dlog [b]$ so $f(Va)=\wp (Va)\dlog
[b]=(FVa-Va)\dlog [b]=0$. \qed

By Lemma 3.5 we have a surjective morphism $W_n(K)\otimes
W_n(K)\to\Omega^1(W_n(K))$, given by $a\otimes b\mapsto a\diff
b$, and its kernel is $\langle a\otimes bc-ab\otimes c-ac\otimes b\,
:\, a,b,c\in W_n(K)\rangle$. Hence $G_n=\Omega^1(W_n(K))/M_n$ also
writes as $G_n=(W_n(K)\otimes W_n(K))/N_n$, where $N_n\sbq
W_n(K)\otimes W_n(K)$ is the preimage of
$M_n\sbq\Omega^1(W_n(K))$.

Explicitly, $N_n$ is the group generated by
\begin{center}
	$a\otimes bc-ab\otimes c-ac\otimes b$, with $a,b,c\in
        W_n(K),$
	
	$Fa\otimes b-a\otimes Vb$, with $a,b\in W_n(K)$,\quad
        $\wp ([a])[b]^{-1}\otimes [b]$, with $a\in K$, $b\in\kk$.
\end{center}

If we write $G_n$ as $(W_n(K)\otimes W_n(K))/N_n$ then
$\alpha_{p^n}:G_n\to\Brr{p^n}(K)$ is given by $a\otimes b\mapsto
((a,b))_{p^n}$.

Similarly, $G'_n=\Omega^1(W_n(K))/M'_n$ also writes as
$G_n=(W_n(K)\otimes W_n(K))/N'_n$, where $N'_n$ has the same
generators as $N_n$ except $\wp ([a])[b]^{-1}\otimes [b]$, with $a\in
K$, $b\in\kk$.

\blm In $G_n=(W_n(K)\otimes W_n(K))/N_n$ we have:

$Fa\otimes b=a\otimes Vb$,\quad $Va\otimes b=a\otimes Fb$,\quad
$a\otimes b=-b\otimes a$\quad $\forall a,b\in W_n(K)$.

$a\otimes b^k=kab^{k-1}\otimes b$\quad $\forall a,b\in W_n(K)$,
$k\geq 0$.

$\wp (a)[b]^{-1}\otimes [b]=0$\quad $\forall a\in W_n(K)$,
$b\in\kk$.

All relations above, except the last one, also hold in
$G'_n=(W_n(K)\otimes W_n(K))/N'_n$.
\elm
\pf When we regard $G_n$ and $G'_n$ as $\Omega^1(W_n(K))/M_n$ and
$\Omega^1(W_n(K))/M'_n$, the relations we want to prove write as
$Fa\diff b=a\diff Vb$, $Va\diff b=a\diff Fb$, $a\diff b=-b\diff
a$, $a\diff b^k=kab^{k-1}\diff b$ and $\wp (a)[b]^{-1}\diff [b]=\wp
(a)\dlog [b]=0$. They follow from Lemma 4.4 and the properties of the
differentials. (For $a\diff b=-b\diff a$ note that $a\diff b+b\diff
a=\diff (ab)\in M'_n\sbq M_n$.) \qed

Since $W_n(K)=W(K)/V^n(W(K))$ the elements of $W_n(K)$ can be
regarded as classes of elements in $W(K)$. Therefore elements of
$W_n(K)\otimes W_n(K)$ can be regarded as classes of elements
of $W(K)\otimes W(K)$. Hence we may regard $\alpha_{p^n}$ as being
defined on $W(K)\otimes W(K)$.

\blm With the convention above, if $m\geq n$ then
$\alpha_{p^n}=p^{m-n}\alpha_{p^m}$, in the sense that
$\alpha_{p^n}(\eta )=p^{m-n}\alpha_{p^m}(\eta )$ $\forall\eta\in
W(K)\otimes W(K)$.
\elm
\pf Since $\alpha_{p^m}$ and $\alpha_{p^n}$ are linear it is enough to
verify our statement when $\eta$ is a generator of $W(K)\otimes W(K)$,
$\eta =a\otimes b$, with $a,b\in W(K)$. But by Corollary 3.9 we have
$((a,b))_{p^n}=p^{m-n}((a,b))_{p^m}$, i.e. $\alpha_{p^n}(a\otimes
b)=p^{m-n}\alpha_{p^m}(a\otimes b)$. \qed

\blm If $n\geq 2$ then $V\otimes V:W_{n-1}(K)\otimes W_{n-1}(K)\to
W_n(K)\otimes W_n(K)$ induces a linear map $f_n:G_{n-1}\to G_n$.

We have $f_n(\eta)=p\eta$ for any $\eta\in W(K)\otimes W(K)$ so $\Ima
f_n=pG_n$. 
\elm
\pf We have $G_{n-1}=(W_{n-1}(K)\otimes W_{n-1}(K))/N_{n-1}$ and
$G_n=(W_n(K)\otimes W_n(K))/N_n$ so we must prove that
$(V\otimes V)(N_{n-1})\sbq N_n$. Equivalently, if $h$ is the
composition $W_{n-1}(K)\otimes W_{n-1}(K)\xrightarrow{V\otimes V}
W_n(K)\otimes W_n(K)\to (W_n(K)\otimes W_n(K))/N_n=G_n$, then we must
prove that $N_{n-1}\sbq\ker h$. 

First note that for any $a,b\in W(K)$, by Lemma 4.5, in $G_n$ we have
$h(a\otimes b)=Va\otimes Vb=FVa\otimes b=pa\otimes b$. More
generally, by linearity we have $h(\eta )=p\eta$ $\forall\eta\in
W(K)\otimes W(K)$. 

To prove that $N_{n-1}\sbq\ker h$ we note that $N_{n-1}$ is generated
by the elements $\eta =a\otimes bc-ab\otimes c-ac\otimes b$, $Fa\otimes
b-a\otimes Vb$ and $\wp([a])[b]^{-1}\otimes [b]$ of
$W(K)\otimes W(K)$. In each case we have $\eta =0$ in
$G_n=(W_n(K)\otimes W_n(K))/N_n$ so $h(\eta )=p\eta =0$, as
well. Hence $N_{n-1}\sbq\ker h$.

We obtain a linear map $f_n:G_{n-1}=(W_{n-1}(K)\otimes
W_{n-1}(K))/N_{n-1}\to G_n=(W_n(K)\otimes W_n(K))/N_n$ given by
$f_n(\eta )=h(\eta )=p\eta$ $\forall\eta\in W(K)\otimes W(K)$. \qed

\blm Let $f:A\to B$ be a morphism of abelian groups. If $A'\sbq A$ is
a subgroup and $B'=f(A')$ and $f'':A/A'\to B/B'$ is the morphism given
by $f''(x+A')=f(x)+B'$ then $\ker f''$ is the image in $A/A'$ of $\ker f\sbq A$.
\elm
\pf We use the snake lemma for the following diagram. 
$$\begin{array}{ccccccc}
0\to & A'&\to & A & \to & A/A' & \to 0\\
{} & \quad\downarrow f' & {} & \quad\downarrow f & {} &
\quad\downarrow f'' & {}\\
0\to & B'&\to & B & \to & B/B' & \to 0
\end{array},$$
where $f'$ is the restriction of $f$ to $A'$. We have $f(A')=B'$ so
$f'$ is surjective so $\coker f'=0$. Then, as a part of the long exact
sequence obtained by the snake lemma, we have the exact sequence $\ker
f\to\ker f''\to\coker f'=0$. Hence $\ker f''$ is the image in $A/A'$ of $\ker
f\sbq A$, as claimed. \qed

\blm If $n\geq 1$ and $T:W_n(K)\to W_1(K)=K$ is the truncation map
then $T\otimes T:W_n(K)\otimes W_n(K)\to K\otimes K$ induces a
surjective morphism $g_n:G_n\to G_1$ satisfying $\ker g_n=pG_n$.
\elm
\pf Since $T$ is given by $(a_0,a_1,\ldots )\mapsto a_0$ its kernel
consists of elements of the form $(0,a_1,a_2,\ldots )$, i.e. $\ker
T=\Ima V$. The truncation map $T\otimes T$ sends the generators of
$N_{n-1}$ to the generators of $N_1$ so $(T\otimes T)(N_n)=N_1$. Let
$g_n:(W_n(K)\otimes W_n(K))/N_n=G_n\to (K\otimes K)/N_1=G_1$
be the induced morphism. Since $T\otimes T$ is surjective, so is
$g_n$. By Lemma 4.8, $\ker g_n$ is the image of $\ker (T\otimes T)\sbq
W_n(K)\otimes W_n(K)$ in $G_n=(W_n(K)\otimes W_n(K))/N_n$. But $\ker 
(T\otimes T)=\ker T\otimes W_n(K)+W_n(K)\otimes\ker T=\Ima V\otimes
W_n(K)+W_n(K)\otimes\Ima V$. 

We now prove that $\ker g_n=pG_n$. Now $pG_n$ is generated by elements
of the form $pa\otimes b$. But $pa\otimes b=V(Fa)\otimes b\in\ker
g_n$. Conversely, we must prove that the generators $Va\otimes b$ and
$a\otimes Vb$ of $\ker g_n$ belong to $pG_n$. If
$b=(b_0,b_1,b_2,\ldots )$ then $b=[b_0]+Vb'$, where
$b'=(b_1,b_2,\ldots )$. Hence $Va\otimes b=Va\otimes [b_0]+Va\otimes
Vb'$. But, by Lemma 4.5, in $G_n$ we have $Va\otimes [b_0]=a\otimes
F[b_0]=a\otimes [b_0]^p=pa[b_0]^{p-1}\otimes [b_0]\in pG_n$ and
$Va\otimes Vb'=FVa\otimes b'=pa\otimes b'\in pG_n$. So $Va\otimes
b\in pG_n$. Similarly, $Vb\otimes a\in pG_n$ so $a\otimes
Vb=-Vb\otimes a\in pG_n$. \qed

\btm The map $\alpha_{p^n}:G_n\to\Brr{p^n}(K)$ is an isomorphism.
\etm
\pf We use the induction on $n$. The case $n=1$ was handled by Lemma
4.3.

Suppose now that $n\geq 2$. By Lemma 4.9, $g_n$ is surjective with
$\ker g_n=pG_n$. By Lemma 4.7, $\Ima f_n=pG_n$. So we have the exact
sequence $G_{n-1}\xrightarrow{f_n}G_n\xrightarrow{g_n}G_1\to 0$. We
also have the obvious exact sequence
$0\to\Brr{p^{n-1}}(K)\to\Brr{p^n}(K)\xrightarrow{p^{n-1}}\Brr
p(K)$. By Lemma 4.6, for any $\eta\in W(K)\otimes W(K)$ we have
$\alpha_{p^{n-1}}(\eta )=p\alpha_{p^n}(\eta )$ and $\alpha_p(\eta
)=p^{n-1}\alpha_{p^n}(\eta )$. Also by Lemma 4.7 $f_n(\eta )=p\eta$,
while $g_n$ is given by truncations so $g_n(\eta )=\eta$. Hence
$\alpha_p(g_n(\eta ))=p^{n-1}\alpha_{p^n}(\eta )$ and
$\alpha_{p^n}(f_n(\eta ))=p\alpha_{p^n}(\eta
)=\alpha_{p^{n-1}}(\eta )$. So we have the commutative exact diagram
$$\begin{array}{ccccccc}
{}& G_{n-1}& \xrightarrow{f_n}& G_n& \xrightarrow{g_n}& G_1&\to 0\\
{}& \quad\downarrow\alpha_{p^{n-1}}& {}& \quad\downarrow\alpha_{p^n}&
{}& \quad\downarrow\alpha_{p}& {}\\
0\to & \Brr{p^{n-1}}(K)& \to & \Brr{p^n}(K) & \xrightarrow{p^{n-1}} &
\Brr p(K) & {}
\end{array}.$$
By the snake lemma we have the exact sequences
$\ker\alpha_{p^{n-1}}\to\ker\alpha_{p^n}\to\ker\alpha_p$ and
$\coker\alpha_{p^{n-1}}\to\coker\alpha_{p^n}\to\coker\alpha_p$. But by
the induction hypothesis $\alpha_{p^{n-1}}$ and $\alpha_p$ are
isomorphisms so their $\ker$ and $\coker$ are zero. It follows that
$\ker\alpha_{p^n}=\coker\alpha_{p^n}=0$ so $\alpha_{p^n}$ is an
isomorphism as well. \qed

As seen from the proof of Theorem 4.10, we have the commutative square
$$\begin{array}{ccc}
G_{n-1}& \xrightarrow{\alpha_{p^{n-1}}}& \Brr{p^{n-1}}(K)\\
\quad\downarrow f_n& {} &\downarrow\\
G_n& \xrightarrow{\alpha_{p^n}}& \Brr{p^n}(K)
\end{array}.$$
Then we have an isomorphism
$\alpha_{p^\infty}:G_\infty\to\Brr{p^\infty}(K)$, where the
$G_\infty=\varinjlim G_n$. We have $G_n=(W_n(K)\otimes
W_n(K))/N_n$ so $G_\infty =\varinjlim (W_n(K)\otimes
W_n(K))/\varinjlim N_n$. Recall that $f_n:G_{n-1}\to G_n$ is induced
by $V\otimes V:W_{n-1}(K)\otimes W_{n-1}(K)\to W_n(K)\otimes
W_n(K)$. But the limit of the directed system $W_1(K)\xrightarrow
VW_2(K)\xrightarrow VW_3(K)\xrightarrow V\cdots$ is $CW(K)$, with the
canonical morphisms $\psi_n:W_n(K)\to CW(K)$ given by
$(a_0,\ldots,a_{n-1})\mapsto (\ldots,0,0,a_0,\ldots,a_{n-1})$. Hence
$\varinjlim W_n(K)\otimes W_n(K)=CW(K)\otimes CW(K)$, with the canonic
morphisms $\psi_n\otimes\psi_n:W_n(K)\otimes W_n(K)\to CW(K)\otimes
CW(K)$. Since $\alpha_{p^n}$ is given by $a\otimes b\mapsto
((a,b))_{p^n}$ $\forall a,b\in W_n(K)$, $\alpha_{p^\infty}$ will be
given by $a\otimes b\mapsto ((a,b))_{p^\infty}$ $\forall a,b\in
CW(K)$, where $((\cdot,\cdot ))_{p^\infty}$ is defined as in
Proposition 3.17.

By the same proof from Lemma 4.7, but with the part involving $\wp
([a])[b]^{-1}\otimes [b]$ ignored, we have $(V\otimes V)(N'_{n-1})\sbq
N'_n$ so $V\otimes V:W_{n-1}(K)\otimes W_{n-1}(K)\to W_n(K)\otimes
W_n(K)$ induces a morphism $f'_n:G'_{n-1}\to G'_n$. So the groups
$G'_n$ also form a directed system and we denote $G'_\infty
=\varinjlim G'_n$. Same as for $G_\infty$, we have $G'_\infty
=CW(K)\otimes CW(K)/\varinjlim N'_n$.

\blm We have $\varinjlim N_n=\varinjlim N'_n$ so $G_\infty
=G'_\infty$.
\elm
\pf We have $\varinjlim N_n=\bigcup_{n\geq
1}(\psi_n\otimes\psi_n)(N_n)$ and $\varinjlim
N'_n=\bigcup_{n\geq 1}(\psi_n\otimes\psi_n)(N'_n)$. Then
$\varinjlim N_n\spq\varinjlim N'_n$ follows from $N_n\spq N'_n$
$\forall n$. Conversely, we must prove that
$(\psi_n\otimes\psi_n)(N_n)\sbq\varinjlim N'_n$ for any $n\geq
1$. It suffices to prove that $(\psi_n\otimes\psi_n)(\eta
)\in\varinjlim N'_n$ for every generator $\eta$ of $N_n$. If $\eta
=a\otimes bc-ab\otimes c-ac\otimes b$ for some $a,b,c\in W_n(K)$
or $Fa\otimes b-a\otimes Vb$ for some $a,b\in W_n(K)$ then
$\eta\in N'_n$ so $(\psi_n\otimes\psi_n)(\eta
)\in(\psi_n\otimes\psi_n)(N'_n)\sbq\varinjlim N'_n$. Assume
now that $\eta =\wp ([a])[b]^{-1}\otimes [b]$ for some $a,b\in K$,
$b\neq 0$. We have $(\psi_n\otimes\psi_n)(\eta
)=(\psi_{n+1}\otimes\psi_{n+1})((V\otimes V)(\eta ))$. We prove that
$(V\otimes V)(\eta)\in N'_{n+1}$ so
$(\psi_n\otimes\psi_n)(\eta
)\in(\psi_{n+1}\otimes\psi_{n+1})(N'_{n+1})\sbq\varinjlim N'_n$. Now
$\eta =[a^pb^{-1}]\otimes [b]-[ab^{-1}]\otimes [b]$ so $(V\otimes
V)(\eta )=V[a^pb^{-1}]\otimes V[b]-V[ab^{-1}]\otimes V[b]$. But by
Lemma 4.5 in $G'_{n+1}$ we have $V[a^pb^{-1}]\otimes
V[b]=FV[a^pb^{-1}]\otimes [b]=p[a^pb^{-1}]\otimes [b]$ and
$V[ab^{-1}]\otimes V[b]=F[ab^{-1}]\otimes F[b]=[a^pb^{-p}]\otimes
[b]^p=p[a^pb^{-p}][b]^{p-1}\otimes [b]=p[a^pb^{-1}]\otimes [b]$. Thus
$V[a^pb^{-1}]\otimes V[b]=V[ab^{-1}]\otimes V[b]$ in
$G'_{n+1}=(W_{n+1}(K)\otimes W_{n+1}(K))/N'_{n+1}$ so $(V\otimes
V)(\eta )=V[a^pb^{-1}]\otimes V[b]-V[ab^{-1}]\otimes V[b]\in
N'_{n+1}$, as claimed. \qed

\blm We have $\varinjlim N_n=N$, where $N\sbq CW(K)\otimes CW(K)$ is
generated by
\begin{center}
	$a\otimes b+b\otimes a$ and $Fa\otimes b-a\otimes Vb$,
        with $a,b\in CW(K)$
	
	$[a]_l\otimes [bc]_l+[b]_l\otimes [ac]_l+[c]_l\otimes [ab]_l$,
        with $a,b,c\in K$, $l\leq 0$. 
\end{center}
Here for $l\leq 0$ by $[a]_l$ we mean $(\ldots,0,0,a,0,\ldots,0)\in
CW(K)$, with $a$ on the $l$th position. Alternatively,
$[a]_{-n}=\psi_{n+1}([a])$ for $n\geq 0$.
\elm
\pf We first prove that $N\sbq\varinjlim N_n$. We must prove that
every generator $\eta$ of $N$ belongs to  $\varinjlim N_n$. If $\eta
=a\otimes b+b\otimes a$ or $Fa\otimes b-a\otimes Vb$ for some
$a,b\in CW(K)$ then for $n\geq 1$ large enough we have
$a=\psi_n(c)$, $b=\psi_n(d)$ for some $c,d\in W_n(K)$. It
follows that $\eta =(\psi_n\otimes\psi_n)(\nu )$, where $\nu
=c\otimes d+d\otimes c$ or $Fc\otimes d-c\otimes Vd$,
respectively. But, by Lemma 4.5, in both cases we have $\nu=0$ in
$G_n$ so $\nu\in N_n$ and $\eta\in
(\psi_n\otimes\psi_n)(N_n)\sbq\varinjlim N_n$. If $\eta
=[a]_l\otimes [bc]_l+[b]_l\otimes [ac]_l+[c]_l\otimes [ab]_l$, for
some $a,b,c\in K$ and $l\leq 0$ then $l=-(n-1)$ for some $n\geq
1$. We have $[\alpha ]_l=\psi_n([\alpha ])$ $\forall\alpha\in
K$. Hence  $\eta =(\psi_n\otimes\psi_n)(\nu )$, where $\nu
=[a]\otimes [bc]+[b]\otimes [ac]+[c]\otimes [ab]$. But, by Lemma 4.5,
in $G_n$ we have $\nu =[a]\otimes [bc]-[ab]\otimes [c]-[ac]\otimes
[b]=[a]\otimes [b][c]-[a][b]\otimes [c]-[a][c]\otimes [b]=0$. So again
$\nu\in N_n$ and $\eta\in\varinjlim N_n$. 

Before proving the reverse inclusion, note that in $(CW(K)\otimes
CW(K))/N$ we have $a\otimes b=-b\otimes a$ and $Fa\otimes
b=a\otimes Vb$ $\forall a,b\in CW(K)$, but also $Va\otimes b=-b\otimes
Va=-Fb\otimes a=a\otimes Fb$.

By Lemma 4.11, we have $\varinjlim N_n=\varinjlim N'_n=\bigcup_{n\geq
  1}(\psi_n\otimes\psi_n)(N'_n)$ so we must prove that
$(\psi_n\otimes\psi_n)(N'_n)\sbq N$ $\forall n\geq 1$. It
suffices to prove that $(\psi_n\otimes\psi_n)(\eta )\in N$ for
all generators $\eta$ of $N_n'$. If $\eta =Fa\otimes b-a\otimes
Vb$ for some $a,b\in W_n(K)$ then $(\psi_n\otimes\psi_n)(\eta
)=Fa'\otimes b'-a'\otimes Vb'\in N$, where $a'=\psi_n(a)$,
$b'=\psi_n(b)$. For elements $\eta$ of the form $a\otimes bc-ab\otimes
c-ac\otimes b$ with $a,b,c\in W_n(K)$ note that the map
$(a,b,c)\mapsto a\otimes bc-ab\otimes c-ac\otimes b$ is trilinear so
it is enough to take the case when $a,b,c$ belong to the set of
generators $\{ V^i[\alpha]\, :\, \alpha\in K,\, 0\leq i\leq n-1\}$ of
$W_n(K)$. So we must prove that $(\psi_n\otimes\psi_n)(\eta )\in N$ for $\eta\in W_n(K)$ of the form 
$$\eta =V^i[a]\otimes V^j[b]V^k[c]-V^i[a]V^j[b]\otimes
V^k[c]-V^i[a]V^k[c]\otimes V^j[b]$$
for some $a,b,c\in K$, $0\leq i,j,k\leq n-1$. But
$V^j[b]V^k[c]=V^{j+k}(F^k[b]F^j[c])=V^{j+k}[b^{p^k}c^{p^j}]$ and
similarly for the other products. Hence
$$\eta =V^i[a]\otimes
V^{j+k}[b^{p^k}c^{p^j}]-V^{i+j}[a^{p^j}b^{p^i}]\otimes
V^k[c]-V^{i+k}[a^{p^k}c^{p^i}]\otimes V^j[b].$$

But for any $\alpha\in K$ we have
$\psi_n([\alpha ])=[\alpha ]_l$, where $l=-(n-1)$, so in
$(CW(K)\otimes CW(K))/N$ we have
\begin{multline*}
(\psi_n\otimes\psi_n)(\eta )\\
=V^i[a]_l\otimes
V^{j+k}[b^{p^k}c^{p^j}]_l-V^{i+j}[a^{p^j}b^{p^i}]_l\otimes
V^k[c]_l-V^{i+k}[a^{p^k}c^{p^i}]_l\otimes V^j[b]_l\\
=V^i[a]_l\otimes V^{j+k}[b^{p^k}c^{p^j}]_l+V^j[b]_l\otimes
V^{i+k}[a^{p^k}c^{p^i}]_l+V^k[c]_l\otimes V^{i+j}[a^{p^j}b^{p^i}]_l.
\end{multline*}
Recall that in $(CW(K)\otimes CW(K))/N$ we have $Fa\otimes
b=a\otimes Vb$ and $Va\otimes b=a\otimes Fb$ $\forall a,b\in
CW(K)$. Therefore
$$V^i[a]_l\otimes V^{j+k}[b^{p^k}c^{p^j}]_l=
F^{j+k}[a]_l\otimes F^i[b^{p^k}c^{p^j}]_l=[a^{p^{j+k}}]_l\otimes
[b^{p^{i+k}}c^{p^{i+j}}]_l.$$
Similarly for the other two terms. Hence if we denote by
$a'=a^{p^{j+k}}$, $b'=b^{p^{i+k}}$, $c'=c^{p^{i+j}}$ then in
$(CW(K)\otimes CW(K))/N$ we have $(\psi_n\otimes\psi_n)(\eta
)=[a']_l\otimes [b'c']_l+[b']_l\otimes [a'c']_l+[c']_l\otimes
[a'b']_l=0$. Thus $(\psi_n\otimes\psi_n)(\eta )\in N$. \qed

In conclusion we have:

\btm With $N$ defined as in Lemma 4.12, there is an isomorphism
$\alpha_{p^\infty}:(CW(K)\otimes CW(K))/N\to\Brr{p^\infty}(K)$, given
by $a\otimes b\mapsto ((a,b))_{p^\infty}$ $\forall a,b\in CW(K)$.
\etm

If in the set of generators of $N$ we replace $a\otimes b+b\otimes a$,
with $a,b\in CW(K)$, by $a\otimes a$, with $a\in CW(K)$, then we obtain a
new subgroup $\overline N\sbq CW(K)\otimes CW(K)$. Since $a\otimes
b+b\otimes a=(a+b)\otimes (a+b)-a\otimes a-b\otimes b\in\overline N$
we have $N\sbq\overline N$. For the reverse inclusion note that
$2a\otimes a=a\otimes a+a\otimes a\in N$ $\forall a\in CW(K)$. If
$p\neq 2$ then for $s$ large enough we also have $p^sa\otimes a=0$ so
$a\otimes a\in N$. If $p=2$ then let $a'\in CW(K)$ such that
$a=Va'$. (Say, $a'=(a,0)$.) Then in $(CW(K)\otimes CW(K))/N$ we have
$a\otimes a=Va'\otimes Va'=FVa'\otimes a'=2a'\otimes a'=0$ so
again $a\otimes a\in N$.

Hence $N=\overline N$. Note that $(CW(K)\otimes CW(K))/\langle
a\otimes a\, :a\,\in CW(K)\rangle =\Lambda^2(CW(K))$. Therefore
Theorem 4.13 can be written in the following equivalent form.

\bco There is an isomorphism
$\alpha_{p^\infty}:\Lambda^2(CW(K))/N'\to\Brr{p^\infty}(K)$, given by
$a\wedge b\mapsto ((a,b))_{p^\infty}$, where
\begin{multline*}
N'=\langle Fa\wedge b-a\wedge Vb\, :\, a,b\in CW(K),\\
[a]_l\wedge [bc]_l+[b]_l\wedge [ac]_l+[c]_l\wedge [ab]_l\, :\,
a,b,c\in K,\, l\leq 0\rangle.
\end{multline*}
\eco
\bigskip

{\bf Acknowledgement} I want to thank Philippe Gille for his interest
in my work and for his advice. It was he who suggested that I should
try to determine the kernel of the map
$\alpha_{p^n}:\Omega^1(W_n(K))\to\Brr{p^n}(K)$ and so to obtain a
representation theorem for $\Brr{p^n}(K)$ which generalizes [GS,
Theorem 9.2.4]. As it turned out, the properties already known,
$((Fa,b))_{p^n}=((a,Vb))_{p^n}$ and $[\wp (a),b)_{p^n}=0$, which
translate to $Fa\diff b-a\diff Vb$, $\wp (a)\dlog
[b]\in\ker\alpha_{p^n}$, were enough to generate $\ker\alpha_{p^n}$.
\bigskip


{\bf References}

[FV] I.V. Fesenko, S.V. Vostokov, "Local Fields and Their Extensions",
The Second Edition, American Math Society, Translations of Math
Monographs vol 121, (2002).

[GS] P. Gille, T. Szamuely, "Central simple algebras and Galois
cohomology" 2nd edition, Cambridge Studies in Advanced Mathematics
165, Cambridge University Press (2017).

[T] Thomas, Lara, "Ramification groups in Artin-Schreier-Witt
extensions", Journal de Th\'eorie des Nombres de Bordeaux 17 (2005),
689–720.

[W] Witt, Ernest "Zyklische K\"orper un Algebren der Charakteristik
$p$ vom Grad $p^n$", (1934).

\bigskip

Institute of Mathematics Simion Stoilow of the Romanian Academy, Calea
Grivitei 21, RO-010702 Bucharest, Romania.

E-mail address: Constantin.Beli@imar.ro



\end{document}